\documentclass{article}
\usepackage{graphicx} 
\usepackage{amsmath,mathtools,amssymb,bm}
\usepackage[margin=1in]{geometry}
\usepackage{biblatex}
\usepackage{xcolor}
\usepackage{float}
\usepackage[export]{adjustbox}
\usepackage{listings}
\usepackage{multicol}
\usepackage{listings}

\usepackage{hyperref}
\usepackage{authblk}

\hypersetup{
    colorlinks,
    citecolor=black,
    filecolor=black,
    linkcolor=black,
    urlcolor=black
}

\usepackage{tikz}
\usetikzlibrary{positioning}
\usetikzlibrary {shapes.misc}

\usepackage{algorithm}
\usepackage{algpseudocode}

\title{Stochastic generative methods for stable and accurate closure modeling of chaotic dynamical systems}

\author[1]{Emily Williams}
\author[1]{David Darmofal}

\affil[1]{Massachusetts Institute of Technology, MA 02139}

\begin{document}

\maketitle

\begin{abstract}
    Traditional deterministic subgrid-scale (SGS) models are often dissipative and unstable, especially in regions of chaotic and turbulent flow. Ongoing work in climate science and ocean modeling motivates the use of stochastic SGS models for chaotic dynamics. Further, developing stochastic generative models of underlying dynamics is a rapidly expanding field. In this work, we aim to incorporate stochastic integration toward closure modeling for chaotic dynamical systems. Further, we want to explore the potential stabilizing effect that stochastic models could have on linearized chaotic systems. We propose parametric and generative approaches for closure modeling using stochastic differential equations (SDEs). We derive and implement a quadratic diffusion model based on the fluctuations, demonstrating increased accuracy from bridging theoretical models with generative approaches. Results are demonstrated on the Lorenz-63 dynamical system.
\end{abstract}

\section{Introduction}

Discretizations of many multiscale systems representing complex physical phenomena contain too many degrees of freedom to simulate accurately given limited computational resources. A motivating example for this is turbulent flow, which is present in many scientific and engineering applications. An alternative modeling framework to direct numerical simulation (DNS) is large eddy simulation (LES), which has demonstrated capabilities in modeling chaotic flow at a more affordable computational cost \cite{garnier_large_2009, sagaut_large_2006}. LES involves the filtering of the Navier-Stokes equations which leads to the appearance of subgrid-scale (SGS) terms that must be modeled. This SGS model is part of a reduced-order model (ROM) which accounts for the SGS dynamics without explicitly resolving them. LES is based on the idea of scale separation, or filtering \cite{sagaut_large_2006}. Consider a box filter of width $\Delta$ applied to the fine-scale states $\bm{x}(t)$ to give the coarse-scale states,
\begin{align}
    \overline{\bm{x}}(t) = \frac{1}{\Delta} \int_{-\Delta/2}^{\Delta/2} \bm{x}(t + t') \,d t'.
\end{align}
Then, the fluctuations are defined by $\bm{x}' = \bm{x} - \overline{\bm{x}}$. In this work, we consider ordinary differential equations (ODEs) of the form
\begin{align}\label{eq:ode}
    \bm{x}_t = \bm{f}(\bm{x})
\end{align}
where $\bm{x}$ is a $d$-dimensional vector of states. Applying the filter to Eq. \ref{eq:ode} gives
\begin{align}
    \overline{\bm{x}}_t = \bm{f}(\overline{\bm{x}}) + \bm{\tau}(\overline{\bm{x}}, \bm{x})
\end{align}
where the SGS term $\bm{\tau}(\overline{\bm{x}}, \bm{x}) = \overline{\bm{f}(\bm{x})} - \bm{f}(\overline{\bm{x}})$. This expression for the SGS terms assumes we have access to the fine-scale states. In practice, we do not have access to the fine-scale states and therefore require a SGS model for $\bm{\tau}(\overline{\bm{x}},\bm{x}) \approx \bm{\tau}^{\mathrm{SGS}}(\overline{\bm{x}})$. An accurate SGS model is essential because unresolved dynamics can significantly impact the degrees of freedom that are kept. Most traditional SGS models $\bm{\tau}^{\mathrm{SGS}}$ are deterministic functions of $\overline{\bm{x}}$ \cite{smagorinsky_general_1963, lilly_numerical_1962, deardorff_numerical_1970, germano_turbulence_1992, vreman_eddyviscosity_2004}. Developing stable and accurate SGS models for their usage in high-fidelity modeling frameworks in simulating turbulent flow for engineering applications remains an active area of research in the computational fluid dynamics (CFD) community \cite{sagaut_large_2009, moser_statistical_2021, garnier_large_2002, bose_wallmodeled_2018, bodart_wallmodeled_2011, slotnick_cfd_2014}.

\begin{figure}[!htb]
    \centering
    \includegraphics[width=0.9\textwidth]{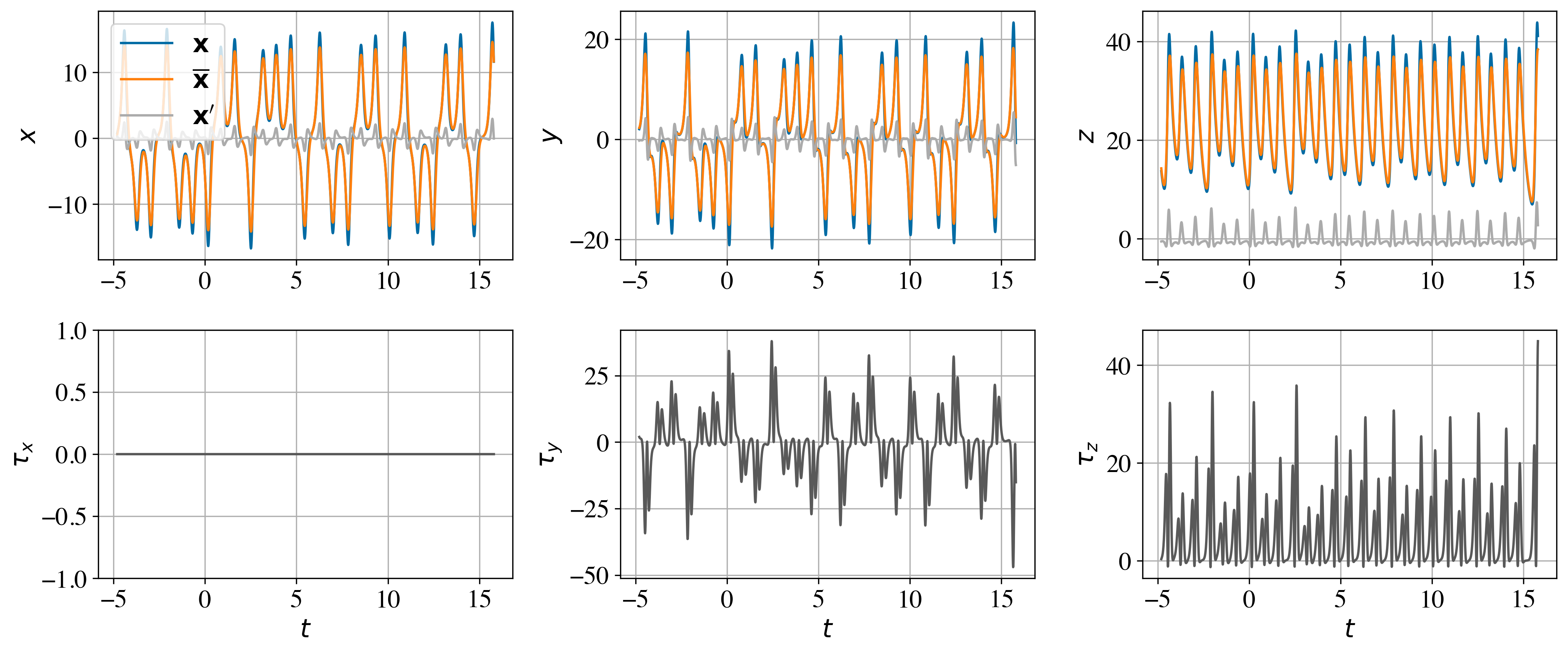}
    \caption{Fine-scale states $\bm{x}$, filtered states $\overline{\bm{x}}$, fluctuations $\bm{x}'$, and subgrid stresses $\bm{\tau}(\overline{\bm{x}},\bm{x})$ for $\Delta = 0.2$.}
    \label{fig:exact_sgs}
\end{figure}

As an example, we will consider the Lorenz-63 chaotic dynamical system defined by
\begin{align}
    \bm{x}_t = \bm{f}(\bm{x}) = \begin{bmatrix}
        \sigma(y-x) \\
        x(r-z) \\
        xy - \beta z
    \end{bmatrix}
\end{align}
with state vector $\bm{x} = (x,y,z)$, initial condition state $\bm{x}(0) = (x_0,y_0,z_0)$ on the attractor, and parameters $\sigma = 10$, $r = 28$, and $\beta = 8/3$. For a filter width of $\Delta = 0.2$ applied to the Lorenz-63 system with states $\bm{x} = (x,y,z)$, the filtered states and exact subgrid states are shown in Figure \ref{fig:exact_sgs}.

\noindent
Stochastic extensions of eddy-viscosity models and explicit algebraic models have since emerged with the motivation of injecting noise to overcome the dissipation present in most deterministic SGS models \cite{adrian_stochastic_1990, leith_stochastic_1990, schumann_stochastic_1997, marstorp_stochastic_2007, rasam_stochastic_2014}. Stochastic processes are used often in mathematical modeling of phenomena that appear to vary chaotically or in a random manner \cite{paul_stochastic_2013, ottinger_stochastic_1996, khasminskii_stochastic_2012}. Stochastic differential equations (SDEs) are ubiquitous in the formulation of these models, including population dynamics, neuron activity, blood clotting, turbulent diffusion, and more. 
Since first introduced in 1976 \cite{hasselmann_stochastic_1976}, there has been ongoing work in climate science and ocean modeling that motivates the use of stochastic SGS models and parameterizations for chaotic dynamics, including stochastic weather and climate models \cite{grudzien_numerical_2020, majda_mathematical_2001, palmer_stochastic_2019, gottwald_stochastic_2017,christensen_simulating_2015, givon_extracting_2004, arnold_stochastic_2013}. While deterministic SGS models can capture the mean response, such models can fail in accounting for the fast-scale variability, leading to an inaccurate representation of the large-scale variability over time \cite{christensen_simulating_2015, grudzien_numerical_2020}. Developing generative models of underlying dynamics from observations is a rapidly expanding field \cite{volkmann2024scalable, finn2024generative, ho2020denoising, song2020score, jacobsen2023cocogen}. Further, probabilistic approaches are being used for predicting subgrid forcing, which can be used to build data-driven stochastic parameterizations \cite{gottwald_datadriven_2016, perezhogin_generative_2023, palmer_stochastic_2019, christensen_simulating_2015, cruzeiro_stochastic_2020}.

\begin{figure}[H]
    \centering
    \includegraphics[width=0.3\textwidth]{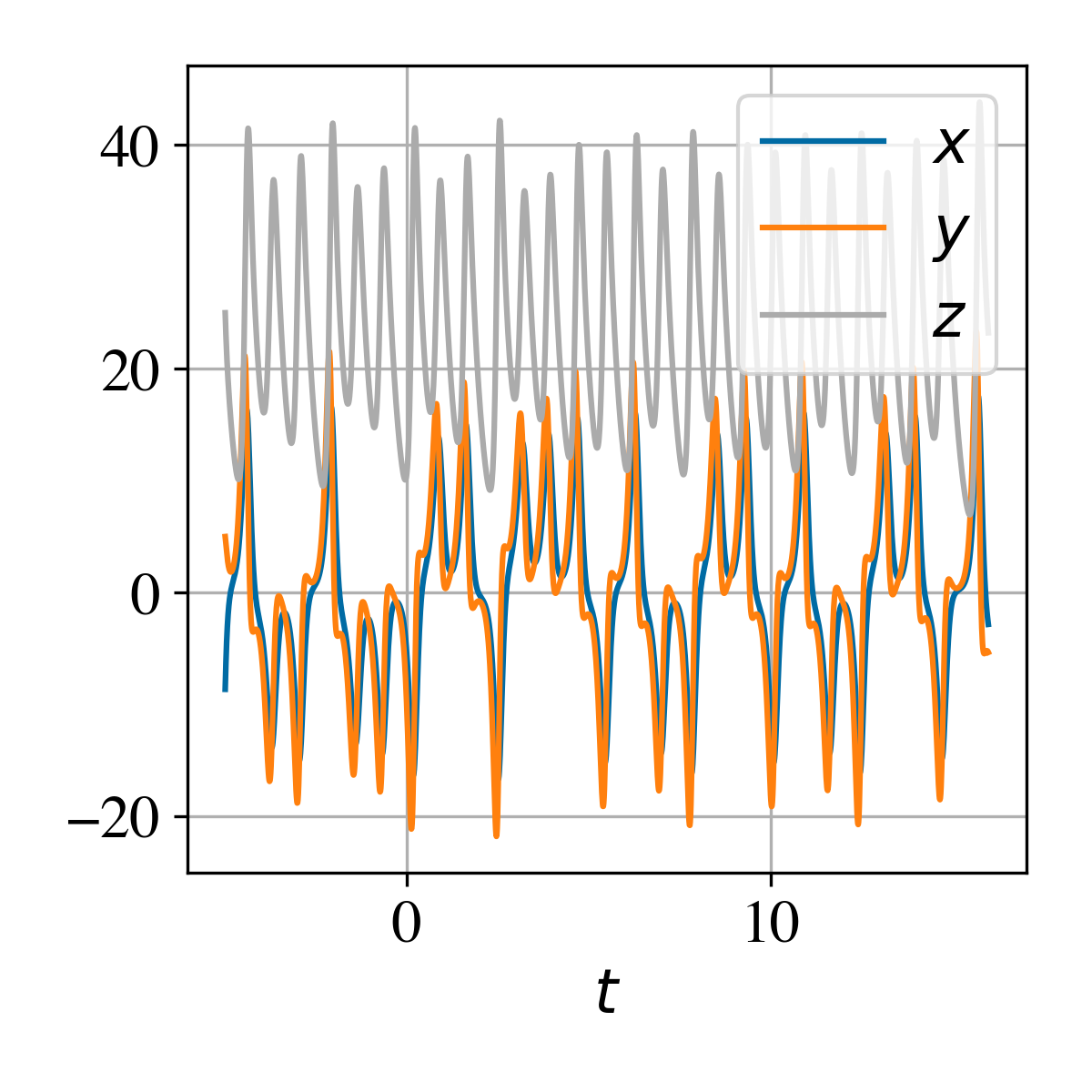}
    \includegraphics[width=0.3\textwidth]{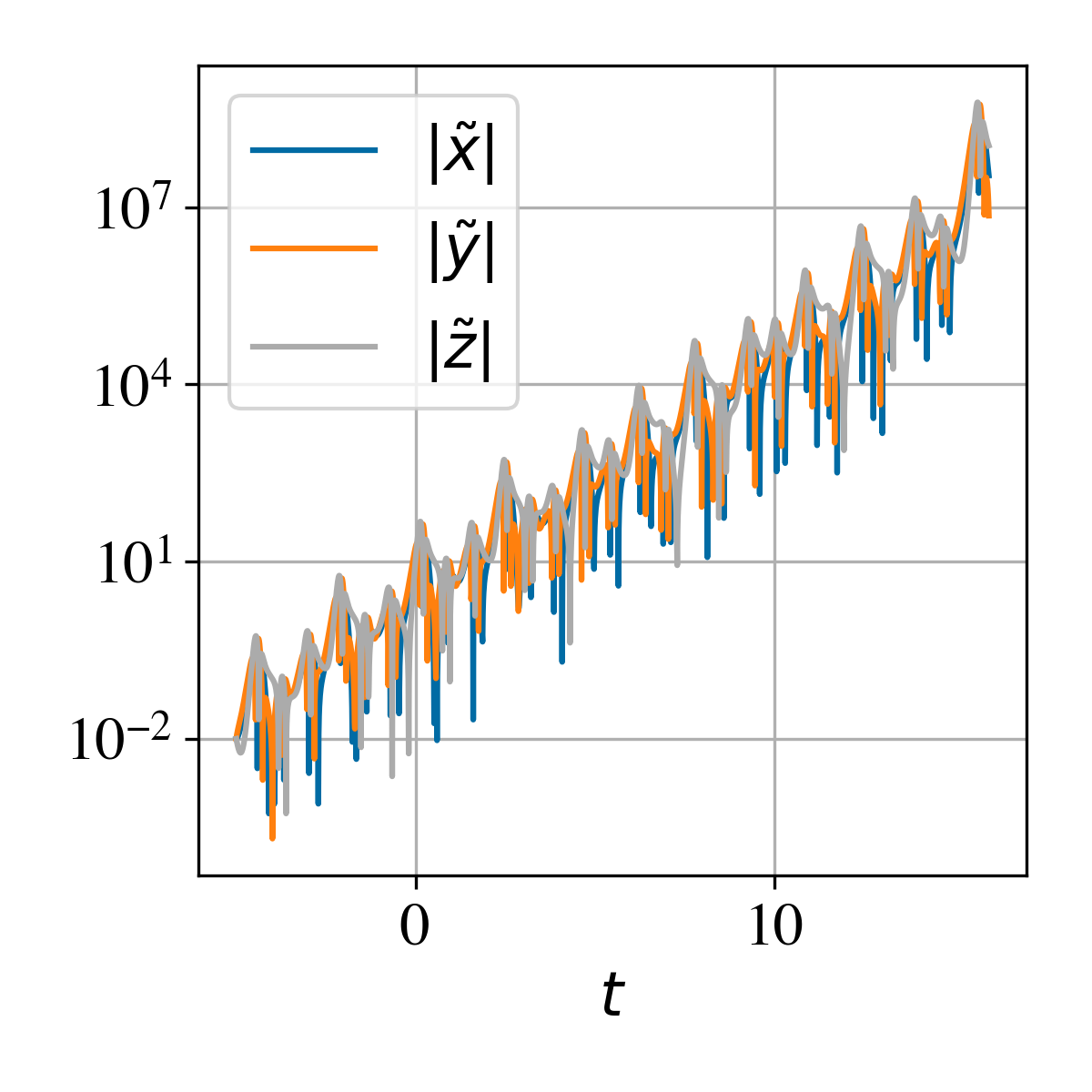}
    \caption{States for (left) nonlinear dynamics and (right) linearized dynamics.}
    \label{fig:lorenz}
\end{figure}

In this work, we first aim to incorporate stochastic integration toward closure modeling for chaotic dynamical systems. Further, we want to explore the potential stabilizing effect that stochastic models could have on the linearized chaotic system,
\begin{align}
    \tilde{\bm{x}}_t = \bm{f_x}(\bm{x}) \tilde{\bm{x}}.
\end{align}
In this case, $\tilde{\bm{x}}$ will experience exponential blow-up.
Figure \ref{fig:lorenz} shows the stable nonlinear states (left) and unstable linearized states (right) for the Lorenz-63 system.

\noindent
Consider now the following SDE for the linearized dynamics,
\begin{align}
    \tilde{\bm{x}}_t = \bm{f_x}(\bm{x}) \tilde{\bm{x}} + \bm{s}(\bm{x},\tilde{\bm{x}},\bm{\xi})
\end{align}
with $\tilde{\bm{x}}(0) = \tilde{\bm{x}}_0$ and stabilization term $\bm{s}(\bm{x},\tilde{\bm{x}},\bm{\xi})$ with stochastic white noise $\bm{\xi} \sim \mathcal{N}(0,1) \in \mathbb{R}^d$. A variety of results exist for stochastic stability \cite{higham_algorithmic_2001, khasminskii_stochastic_2012, liu_stochastic_2019}. Consider a linear SDE 
\begin{align}
    d\tilde{x} = \lambda \tilde{x} dt + \mu \tilde{x} dW
\end{align}
where $dW \sim \mathcal{N}(0,dt)$ is a Brownian increment. This system is asymptotically stable (with probability = 1): $\lim_{t \rightarrow \infty} |\tilde{x}(t)| = 0 \text{ if } \text{Real}\{\lambda - \frac{1}{2}\mu^2 \} < 0$. Thus, for unstable systems with drift $\lambda > 0$, a deterministic system can be stabilized with sufficient diffusion $\mu$. In this work, following the approach from \cite{dietrich_learning_2023}, we formulate probabilistic loss functions and implement parametric stochastic models to demonstrate the applicability of using SDEs for closure and stability.

Motivated by potential dependency of the subgrid terms on the filtered states, we also are interested in using generative approaches for learning probabilistic distributions of the subgrid terms. We formulate and implement guided flow and diffusion models to sample from learned distributions of the subgrid stresses directly. These results demonstrate the benefits of incorporating Langevin dynamics through score matching to better predict the subgrid stresses. Finally, we derive an analytical quadratic model for the subgrid stresses dependent on the fluctuations of the dynamical states. We implement a guided diffusion model conditioned on the filtered states to give the contribution from the fluctuations on the dynamics which are then sampled for the quadratic model for closure. These approaches are demonstrated on the Lorenz-63 chaotic dynamical system. The results contribute to the promising methodology of incorporating stochastic and generative approaches with theoretical models toward improving closure predictions for chaotic systems.

\section{Parametric model}

\subsection{Nonlinear closure}
As a first example, we use the approach and model presented in \cite{dietrich_learning_2023} for the formulation and implementation of the parametric models. Specifically, consider the following SGS model
\begin{align}
    \bm{\tau}(\overline{\bm{x}},\bm{x}) \approx \bm{\tau}^{\mathrm{SGS}} (\overline{\bm{x}},\bm{\xi}) = \bm{\Lambda}(\overline{\bm{x}}) + \bm{\Gamma}(\overline{\bm{x}}) \bm{\xi}
\end{align}
with drift vector $\bm{\Lambda}\in \mathbb{R}^d$, diffusion matrix $\bm{\Gamma}\in \mathbb{R}^{d\times d}$, and white noise $\bm{\xi} \sim \mathcal{N}(0,1) \in \mathbb{R}^d$ characterized by a Wiener process. The Euler-Maruyama \cite{kloeden_numerical_1992, bayram_numerical_2018} integration step is obtained as
\begin{align}
    \overline{\bm{x}}_{n+1} = \overline{\bm{x}}_n + (\bm{f}(\overline{\bm{x}}_n) + \bm{\Lambda}(\overline{\bm{x}}_n)) h + \bm{\Gamma}(\overline{\bm{x}}_n) \Delta \bm{W}_n \label{eq:em_sgs}
\end{align}
where $\Delta\bm{W} \sim \mathcal{N}(0,h) \in \mathbb{R}^d$ is a Brownian increment. For $\overline{\bm{x}}_{n+1}$ conditioned on $(\overline{\bm{x}}_n, h)$ we have
\begin{align}
    \overline{\bm{x}}_{n+1} \sim \mathcal{N}(\overline{\bm{x}}_n + (\bm{f}(\overline{\bm{x}}_n) + \bm{\Lambda}(\overline{\bm{x}}_n)) h, \bm{\Gamma} \bm{\Gamma}^T h) = \mathbb{P}(\overline{\bm{x}}_{n+1} | \overline{\bm{x}}_n, h)
\end{align}
Define $\{\bm{\Lambda}_\theta, \bm{\Gamma}_\theta\}$ as our model with parameters $\theta$ that approximates $\{\bm{\Lambda},\bm{\Gamma}\}$. Following \cite{dietrich_learning_2023}, we define the conditioned probability density of observing $\overline{\bm{x}}_{n+1}$ as $\mathbb{P}_\theta (\overline{\bm{x}}_{n+1} | \overline{\bm{x}}_n, h)$ using the likelihood of the multivariate normal distribution
\begin{align}
    \mathbb{P}_\theta = (2 \pi)^{-d/2} (\det(\bm{M}_\theta))^{-1/2} \exp \left( -\frac{1}{2}(\overline{\bm{x}}_{n+1} - \bm{a}_\theta)^T \bm{M}_\theta^{-1} (\overline{\bm{x}}_{n+1} - \bm{a}_\theta) \right)
\end{align}
where $\bm{a}_\theta = \overline{\bm{x}}_n + (\bm{f}(\overline{\bm{x}}_n) + \bm{\Lambda}_\theta(\overline{\bm{x}}_n)) h$ and $\bm{M}_\theta = \bm{\Gamma}_\theta \bm{\Gamma}_\theta^T h$. We want to maximize the likelihood of having observed the training data generated using the true filtered solution, which is equivalent to minimizing $\mathbb{E}[-\log \mathbb{P}_\theta]$
\begin{align}
    \theta \coloneqq \arg \max_{\hat{\theta}} \mathbb{E}[\log \mathbb{P}_{\hat{\theta}}] = \arg \min_{\hat{\theta}} \mathbb{E}[-\log \mathbb{P}_{\hat{\theta}}] \approx \arg \min_{\hat{\theta}} \frac{1}{N} \sum_{j=1}^N -\log \mathbb{P}_{\hat{\theta}}^j = \arg \min_{\hat{\theta}} \frac{1}{N} \sum_{j=1}^N \mathcal{L}_{\hat{\theta}}^j
\end{align}
The expectation is approximated by the law of large numbers by taking the mean over all $N$ training samples; mathematically, $1/N \sum_{j=1}^N \mathcal{L}_\theta^j$, where $j$ refers to the $j$-th training example. Following \cite{dietrich_learning_2023}, we have the following loss function $\mathcal{L}_\theta$
\begin{align}
    \mathcal{L}(\theta|\overline{\bm{x}}_{n+1},\overline{\bm{x}}_n,h) &= -\log \mathbb{P}_\theta = \frac{d}{2}\log(2\pi) + \frac{1}{2}\log(\det(\bm{M}_\theta)) + \frac{1}{2}(\overline{\bm{x}}_{n+1} - \bm{a}_\theta)^T \bm{M}_\theta^{-1}(\overline{\bm{x}}_{n+1} - \bm{a}_\theta) \label{e:loss_func1}
\end{align} 
To generate the training data, a nominal path is simulated first: $\bm{x}_t = \bm{f}(\bm{x})$ over some time domain, starting from a prescribed initial condition, from $t=0$ to $T$ with timestep $\Delta t$ sufficiently fine for simulating all relevant timescales. Then, the filtered states are computed by approximating the filter integral using a quadrature rule with filter width $\Delta$. We set $h$ and $\Delta$ to be equal, though this is not required. Following Algorithm \ref{alg:train}, we generate $T - K + 1$ training samples. The loss function $\mathcal{L}(\theta|\overline{\bm{x}}_{n+1},\overline{\bm{x}}_n,h)$ and training data $\{\overline{\bm{x}}_i^0, \overline{\bm{x}}_i^h\}_{i=0}^{T-K}$ are used to train $\{\bm{\Lambda}_\theta, \bm{\Gamma}_\theta\}$.

\begin{algorithm}[H]
\caption{Generate training data}\label{alg:train}
\begin{algorithmic}[1]
\State Simulate $\bm{x}_t = \bm{f}(\bm{x})$ from $t=0$ to $T\Delta t$ 
\State Compute $\overline{\bm{x}}^{\mathrm{true}}$ using $\bm{x}$ and $\Delta$ (filter width)
\State $h = \Delta = K \Delta t$
\For{$i = 0,1,\dots,T - K$}
    \State $\overline{\bm{x}}_i^0 = \overline{\bm{x}}^{\mathrm{true}}_i$
    \State $\overline{\bm{x}}_i^h = \overline{\bm{x}}^{\mathrm{true}}_{i+K}$
    \State Save $(\overline{\bm{x}}_i^0,\overline{\bm{x}}_i^h)$, which are associated with $\overline{\bm{x}}^{\mathrm{true}}_i$
\EndFor
\end{algorithmic}
\end{algorithm}
\noindent
The nominal path is simulated with a small $\Delta t = 0.001$ using a second-order discontinuous Galerkin integration scheme to create a set of discrete $\bm{x}$ values representing the path, $\{\bm{x}_n\}_{n=0}^T$. The ResNet from \cite{dietrich_learning_2023}, which is trained to minimize the loss function in Eq. \ref{e:loss_func1}, is used as the model. We assume a diagonal diffusion matrix, i.e., the stochastic noise is uncoupled between the states. The neural network is trained for $\Delta = h = 0.01$. We use a multilayer perceptron (MLP) network with 1 hidden layer with 2 neurons per layer and rectified linear unit (ReLU) activation function.

\begin{figure}[H]
    \centering
    \includegraphics[width=0.9\textwidth]{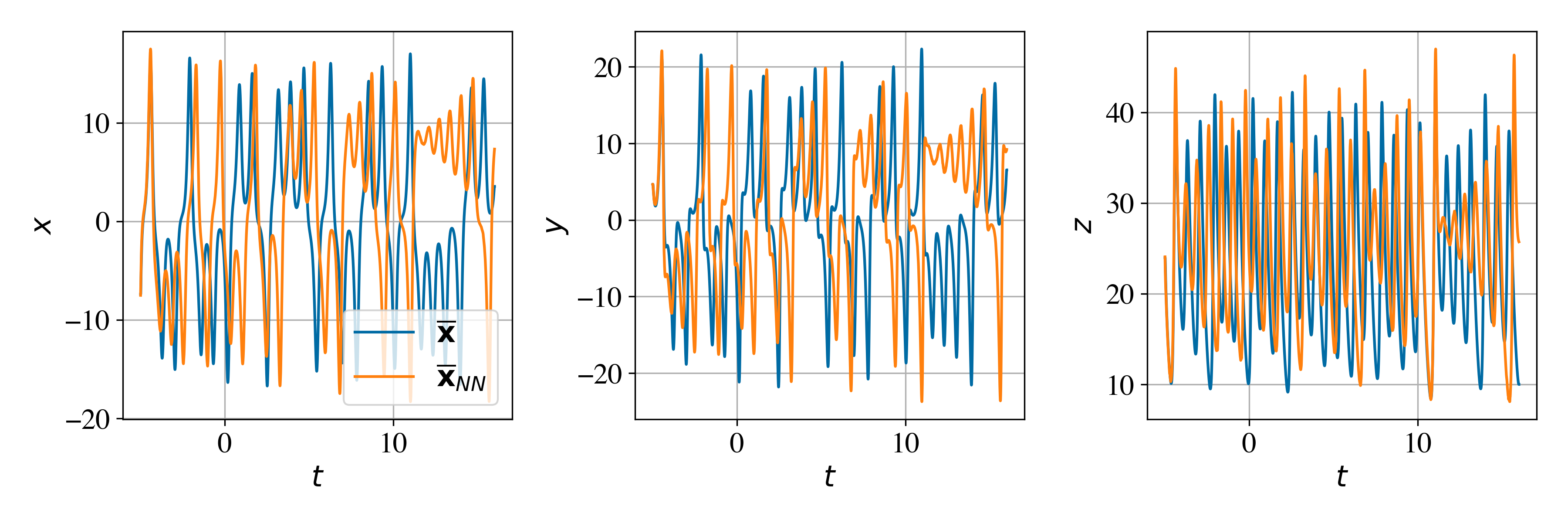}
    \caption{Time evolution of dynamical system with stochastic SGS model.}
    \label{fig:training_NL_states}
\end{figure}

\begin{figure}[!htb]
    \centering
    \includegraphics[width=0.4\textwidth]{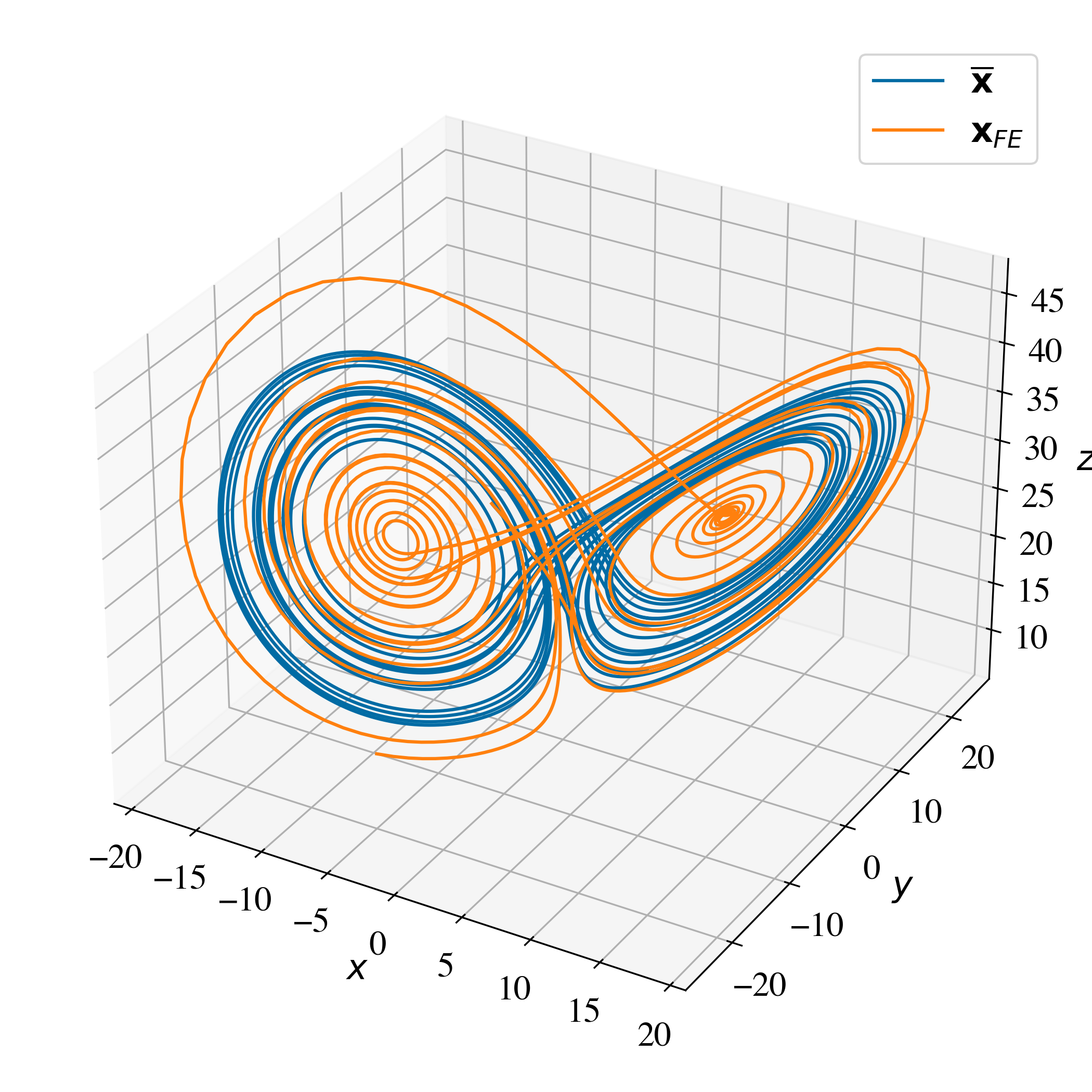}
    \includegraphics[width=0.4\textwidth]{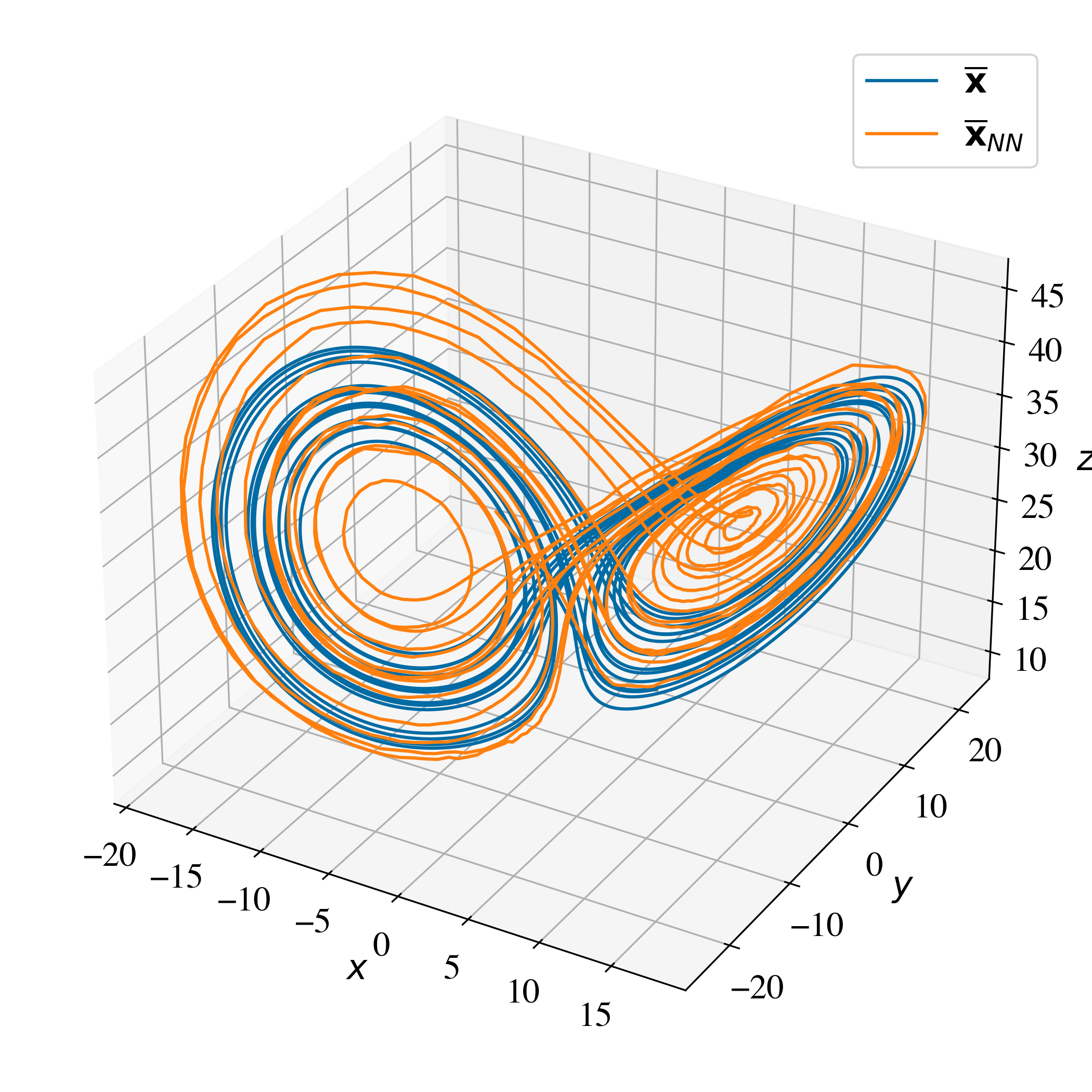}
    \caption{Comparison of (left) deterministic forward Euler $\bm{x}_{FE}$ with timestep $h$ and (right) learned Euler-Maruyama $\overline{\bm{x}}_{NN}$ with timestep $h$ to true filtered states $\overline{\bm{x}}$.}
    \label{fig:training_NL}
\end{figure}

Figure \ref{fig:training_NL_states} shows the evolution of the states with the learned stochastic SGS model from the ResNet model (labeled NN) with timestep $h$ as formulated in Eq. \ref{eq:em_sgs} compared to the true filtered states. Figure \ref{fig:training_NL} (left) shows the trajectory using deterministic forward Euler (FE) integration with timestep $h$. Figure \ref{fig:training_NL} (right) shows the Euler-Maruyama result with the stochastic SGS model. Both trajectories are compared to the true $\overline{\bm{x}}$ from filtering the fine-scale trajectory. The FE and NN trajectories start from the same initial condition as the filtered trajectory. Qualitatively, the trajectory with the stochastic SGS model appears to spend more time near the true filtered trajectory, without stalling around the focal points.

\noindent
Quantitative comparisons can be measured by the Wasserstein distance.  If $P$ and $Q$ are cumulative distribution functions (cdfs), the Wasserstein distance can be defined as
\begin{align}
    W_1 = \int |P(\bm{x}) - Q(\bm{x})| d\bm{x}
\end{align}
A smaller Wasserstein distance means the two distributions are closer to each other, such that the ``cost'' of moving mass from one distribution to the other is minimal \cite{kolouri2017optimalmass,ramdas2015wassersteinsampletestingrelated}. The FE trajectory compared with the true $\overline{\bm{x}}$ has Wasserstein distances of $(x,y,z) = (1.78, 2.20, 3.10)$. The trajectory with the stochastic SGS model compared with the true $\overline{\bm{x}}$ has Wasserstein distances of $(x,y,z) = (1.39, 1.66, 2.02)$. Thus, the trajectory with the stochastic SGS model is closer than the FE trajectory to the true $\overline{\bm{x}}$. Effectively, this demonstrates the feasibility and applicability in using an SDE as a closure model for a chaotic system.

\subsection{Linearized stabilization}
\begin{figure}[H]
    \centering
    \begin{minipage}{.7\textwidth}
        \begin{algorithm}[H]
        \caption{Generate training data}\label{alg:train2}
        \begin{algorithmic}[1]
        \State Simulate $\bm{x}_t = \bm{f}(\bm{x})$ from $t=0$ to $T \Delta t$
        \For{$n = 0,1,\dots,T$}
            \For{$k = 1,\dots,K$}
                \State $\Tilde{\bm{x}}_n^k = \epsilon\,\mathcal{N}(0,1)$ (small random perturbation)
                \State $\bm{x}_n^k = \bm{x}_n + \Tilde{\bm{x}}_n^k$ (nominal value perturbed)
                \State $\bm{x}_{n+1}^k = \bm{x}_n^k + \bm{f}(\bm{x}_n^k)\Delta t$ (one nonlinear step)
                \State $\Tilde{\bm{x}}_{n+1}^k = \bm{x}_{n+1}^k - \bm{x}_{n+1}$ (difference)
                \State Save $(\Tilde{\bm{x}}_n^k,\Tilde{\bm{x}}_{n+1}^k)$, which are associated with $\bm{x}_n$
            \EndFor
        \EndFor
        \end{algorithmic}
        \end{algorithm}
    \end{minipage}%
    \begin{minipage}{0.3\textwidth}
        \includegraphics[width=\textwidth]{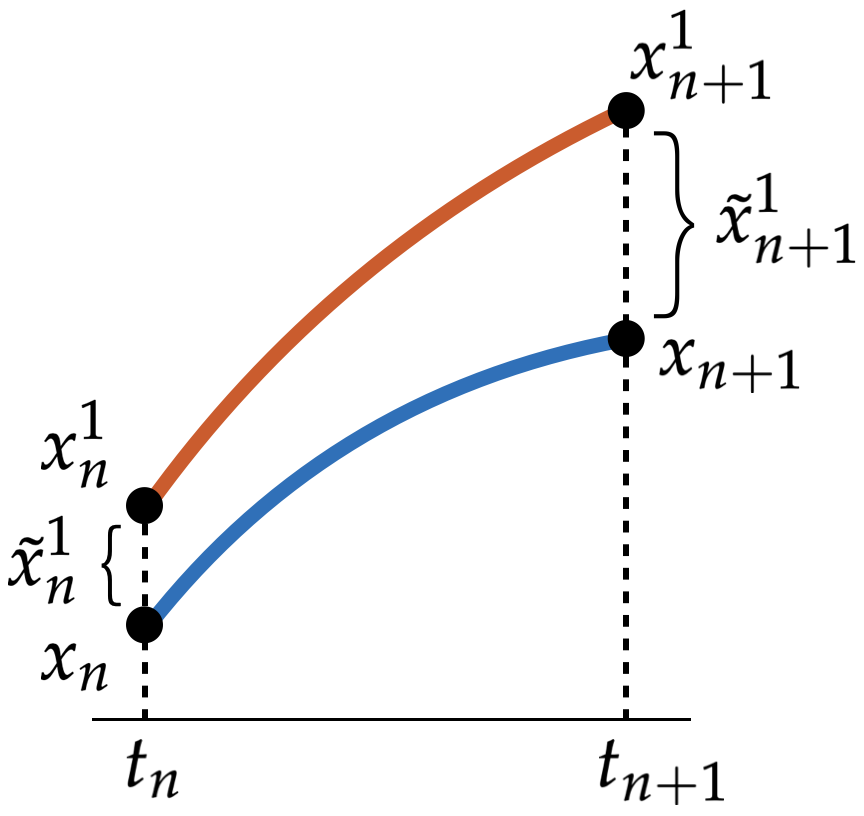}
    \end{minipage}
\end{figure}

We follow a similar approach as when modeling $\bm{\tau}^{\mathrm{SGS}}$ for the nonlinear filtered problem. We seek a stabilization matrix $\bm{\Sigma}$, implicit in the data, to append as the stochastic stabilizing term to the linearized model. The Euler-Maruyama \cite{kloeden_numerical_1992, bayram_numerical_2018} integration step is obtained as
\begin{align}
    \tilde{\bm{x}}_{n+1} = \tilde{\bm{x}}_n + \bm{f_x} (\bm{x}_n) \tilde{\bm{x}}_n \Delta t + \bm{\Sigma}(\bm{x}_n, \tilde{\bm{x}}_n) \Delta \bm{W}_n
\end{align}
where $\bm{\Sigma}$ is a diffusion matrix to stabilize the SDE. For $\Tilde{\bm{x}}_{n+1}$ conditioned on $(\bm{x}_n,\Tilde{\bm{x}}_n,\Delta t)$ we have
\begin{align}
    \Tilde{\bm{x}}_{n+1} \sim \mathcal{N}(\Tilde{\bm{x}}_n + \bm{f_x}(\bm{x}_n)\Tilde{\bm{x}}_n\Delta t,\bm{\Sigma}\bm{\Sigma}^T\Delta t) = \mathbb{P}(\Tilde{\bm{x}}_{n+1} | \bm{x}_n, \Tilde{\bm{x}}_n, \Delta t)
\end{align}
Define $\bm{\Sigma}_\theta$ as our model that approximates $\bm{\Sigma}$ with parameters $\theta$. We define the conditioned probability density of observing $\Tilde{\bm{x}}_{n+1}$, parameterized by the network as $\mathbb{P}_\theta(\Tilde{\bm{x}}_{n+1} | \bm{x}_n, \Tilde{\bm{x}}_n, \Delta t)$. Then, we use the likelihood of the multivariate normal distribution to arrive at the loss function
\begin{align}
    \mathcal{L}_\theta = \frac{d}{2}\log(2\pi) + \frac{1}{2}\log(\det(\bm{M}_\theta)) + \frac{1}{2}(\Tilde{\bm{x}}_{n+1} - \bm{a}_\theta)^T \bm{M}_\theta^{-1}(\Tilde{\bm{x}}_{n+1} - \bm{a}_\theta)
    \label{e:loss_func}
\end{align}
where $\bm{a}_\theta = \Tilde{\bm{x}}_n + \bm{f_x}(\bm{x}_n)\Tilde{\bm{x}}_n\Delta t$ and $\bm{M}_\theta = \bm{\Sigma}_\theta \bm{\Sigma}_\theta^T \Delta t$. A nominal path is simulated to create a set of discrete $\bm{x}$ values representing the path, $\{\bm{x}_n\}_{n=0}^T$. Following Algorithm \ref{alg:train2}, we generate $K \times (T+1)$ training samples, with the labeled schematic. We use the loss function $\mathcal{L}_\theta$ and the training data $\{\{\bm{x}_n,\Tilde{\bm{x}}_n^k,\Tilde{\bm{x}}_{n+1}^k\}_{k=1}^K\}_{n=0}^T$ to train $\bm{\Sigma}_\theta$ for a fixed $\Delta t$. Eq. \ref{e:loss_func} is for a single training point. We want to minimize $\mathbb{E}[\mathcal{L}_\theta]$ across the training set, so we take the mean over all $N = K \times (T+1)$ training samples.

\begin{figure}[H]
    \centering
    \includegraphics[width=0.8\textwidth]{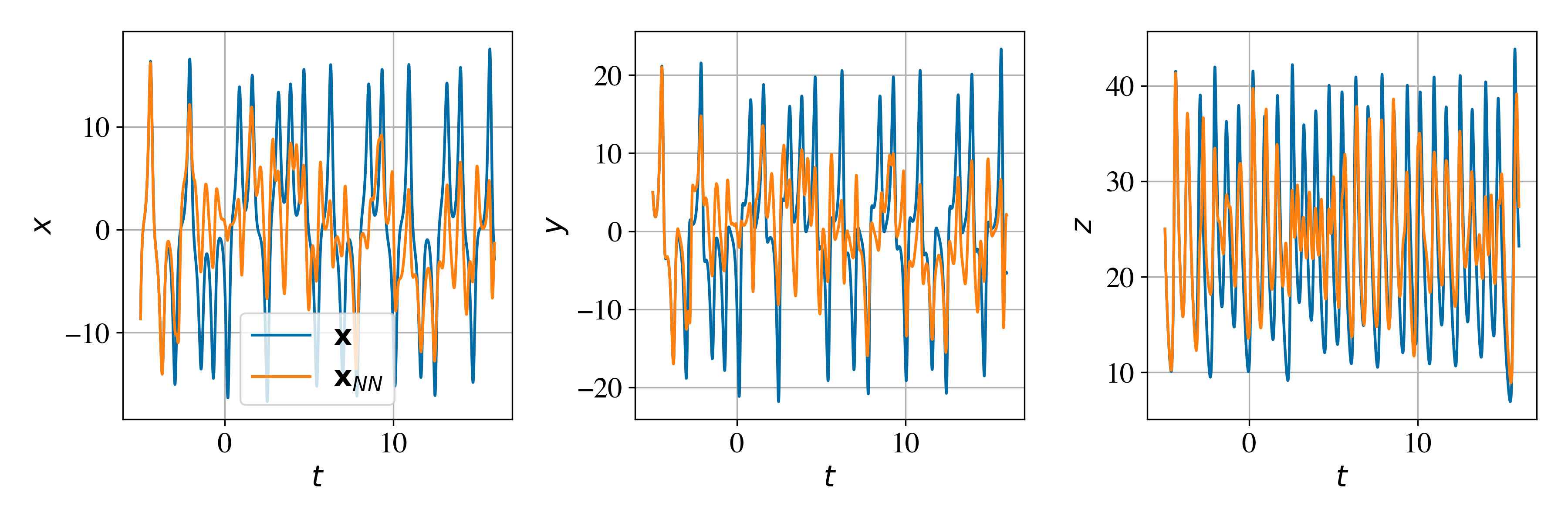}
    \caption{Time evolution of dynamical system with parametric stabilization model.}
    \label{fig:realization}
\end{figure}

\begin{figure}[H]
    \centering
    \includegraphics[width=0.4\textwidth]{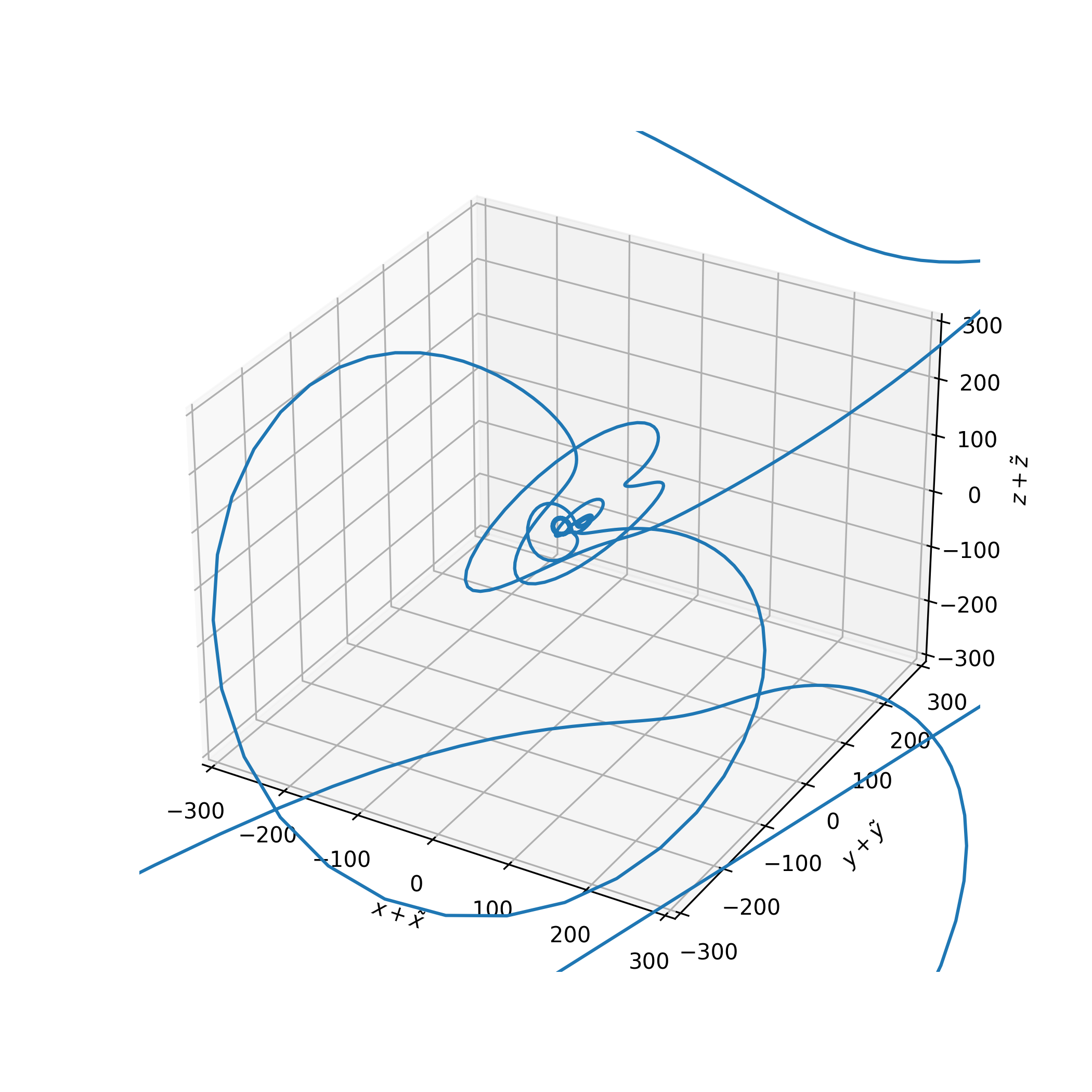}
    \includegraphics[width=0.4\textwidth]{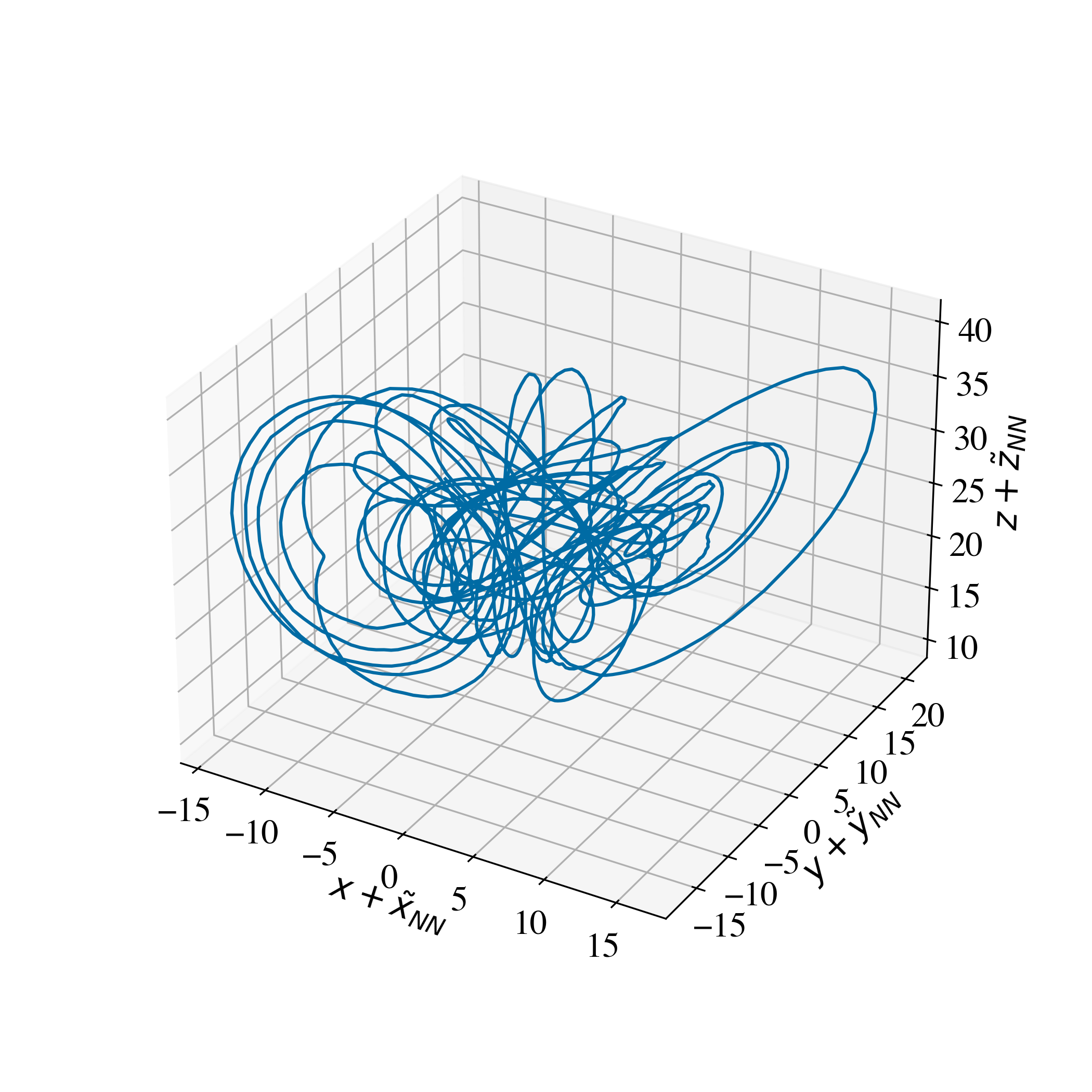}
    \caption{(Left) unstable and (right) stabilized linearized dynamics for parametric stabilization model.}
    \label{fig:training-LIN}
\end{figure}

\noindent
We use the ResNet from \cite{dietrich_learning_2023} and assume that $\bm{\Sigma}$ is a diagonal matrix, i.e., the stochastic noise is uncoupled between the states. The MLP network therefore has 6 inputs and 3 outputs, with 3 hidden layers with 10 neurons per layer and ReLU activation function. The learned parametric model is able to consistently stabilize the dynamics, without injecting excessive diffusion, over a large number of test samples. We can assess the evolution of the states $\bm{x}_{NN} = \bm{x} + \tilde{\bm{x}}_{NN}$. Figure \ref{fig:realization} shows a sample trajectory demonstrating stabilizing transitions between the two focal points, without exponential blow-up. Figure \ref{fig:training-LIN} shows that the parametric model (right) stabilizes the unstable linearized states (left). Effectively, we have demonstrated the feasibility and applicability in using an SDE to stabilize an unstable chaotic system.

\section{Generative model}

\subsection{Nonlinear closure}

A generative approach samples the SGS term from the data distribution $p_{data}$ along the trajectory in time via learning a denoising diffusion process. We assume access to some initial distribution $p_{init}$ that we can easily sample from, such as the Gaussian $p_{init} = \mathcal{N}(0,I_d)$. The goal of our model is then to transform samples from $\bm{t}_0 \in p_{init}$ into samples from $\bm{t}_1 \sim p_{data}$. We assume access to a finite number of examples sampled independently from $p_{data}$, which together serve as a proxy for the true distribution. The desired transformation can be obtained as the simulation of an SDE.

In many cases, we want to generate an object conditioned on some data $\bm{y}$. Guided (conditional) generation involves sampling from $\bm{\tau} \sim p_{data}(\cdot | \bm{y})$, where $\bm{y}$ is a conditioning variable and $p_{data}(\cdot|\bm{y})$ is the conditional data distribution. Whereas for unguided generation we simply want to generate any realistic subgrid stress, we would like to be able to condition on the filtered state. The subgrid stress we would like to generate is denoted as $\bm{\tau} \in \mathbb{R}^3$ and the conditioning variable (in this case, the filtered state) as $\overline{\bm{x}} \in \mathbb{R}^3$. Consider a conditional distribution for the subgrid stresses

\begin{figure}[H]
    \centering
    \includegraphics[width=0.9\textwidth]{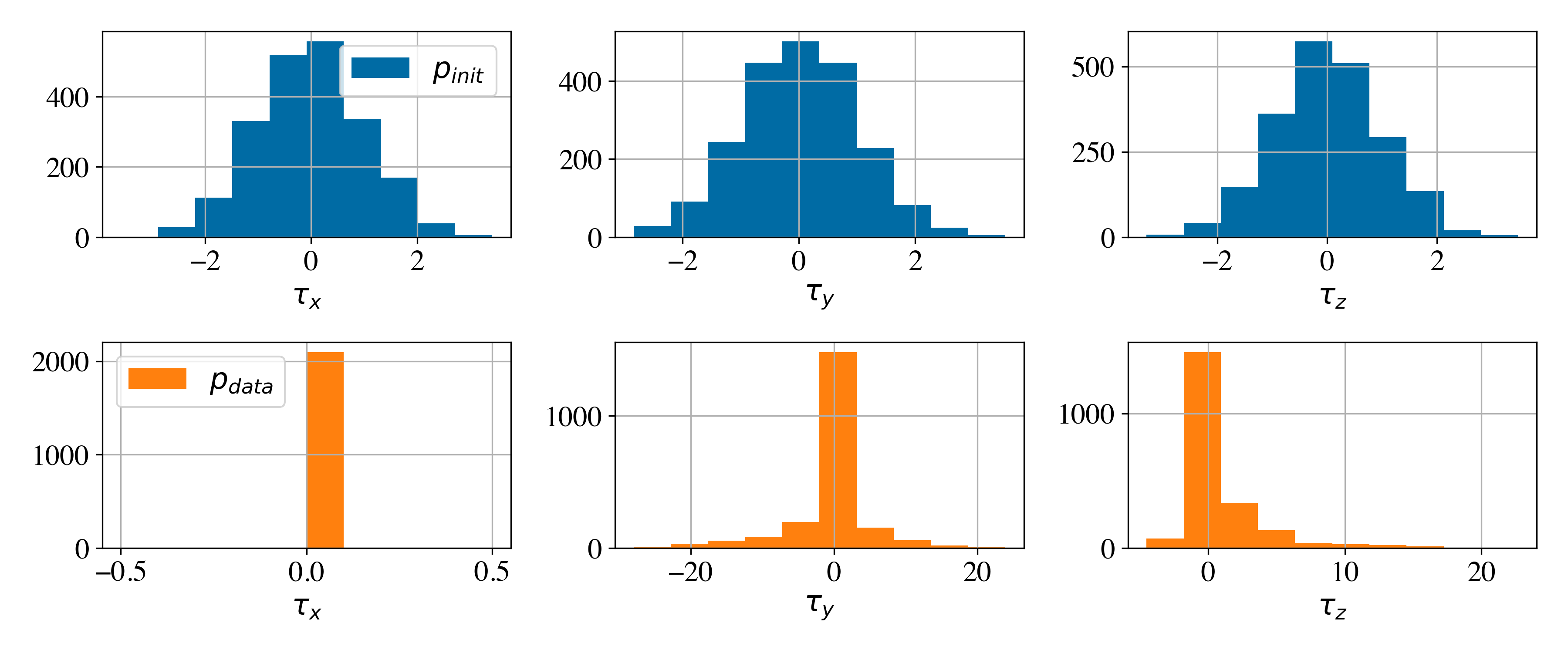}
    \caption{Initial $p_{init}$ and learned $p_{data}$ distributions for guided SGS terms without Langevin dynamics.}
    \label{fig:guided_pdf_flow}
\end{figure}

\begin{figure}[H]
    \centering
    \includegraphics[width=0.9\textwidth]{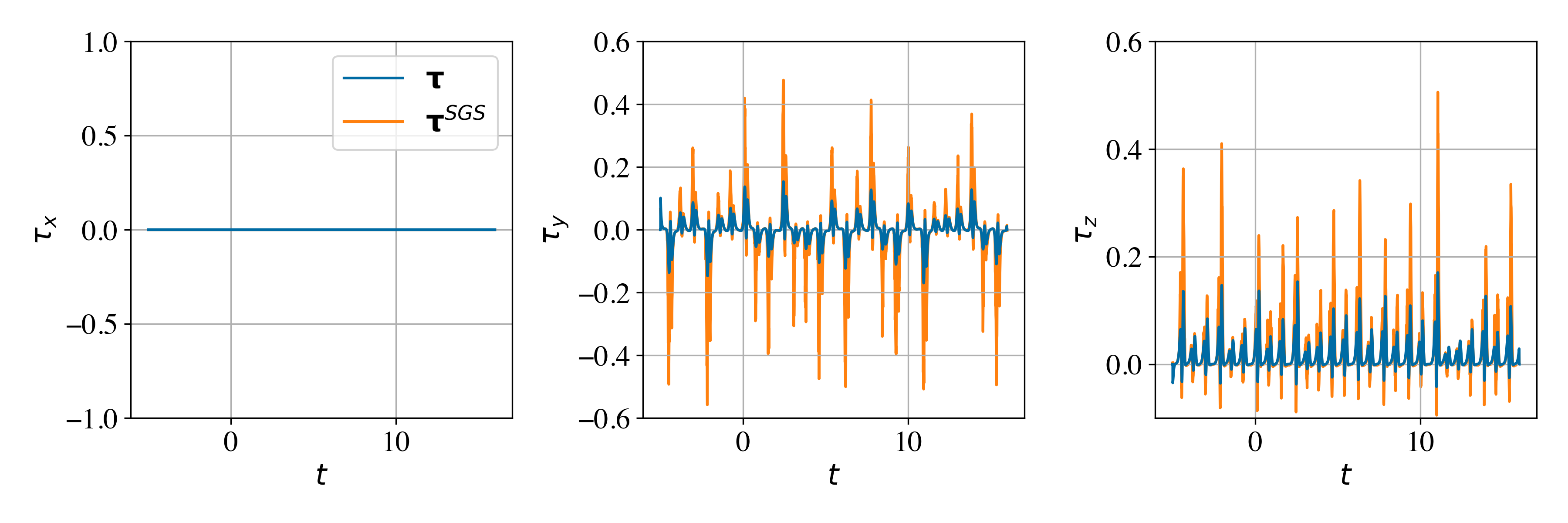}
    \caption{True $\bm{\tau}$ and generated $\bm{\tau}^{\mathrm{SGS}}$ with guided flow model.}
    \label{fig:guided_tau_flow}
\end{figure}

{\fontsize{11}{14}
    \linespread{1.75}
    \begin{center}
    \begin{tikzpicture}[
    squarednode/.style={rectangle, rounded corners, draw=teal!25, fill=cyan!5, minimum height=0.5em, minimum width=1.5em},
    ]

    \path (0,0) node(a) [rectangle, rounded corners, draw=teal!25, fill=cyan!5, minimum height=0.5em, minimum width=1.5em, align=left, inner xsep=0.5em, inner ysep=0.5em]
    {   $\bm{x}_t = \bm{f}(\bm{x})$\\
        $\bm{x}(0) = \bm{x}_0$};
    \path (4.5,0) node(b) [rectangle, rounded corners, draw=teal!25, fill=cyan!5, minimum height=0.5em, minimum width=1.5em, align=left, inner xsep=0.5em, inner ysep=0.5em]
    {   $\overline{\bm{x}}_t = \overline{\bm{f}(\bm{x})}$\\
        $\overline{\bm{x}}(0) = \overline{\bm{x}}_0$};
    \path (10.5,0) node(c) [rectangle, rounded corners, draw=teal!25, fill=cyan!5, minimum height=0.5em, minimum width=1.5em, align=left, inner xsep=0.5em, inner ysep=0.5em] 
    {   $\overline{\bm{x}}_t = \bm{f}(\overline{\bm{x}}) + \bm{\tau}^{\mathrm{SGS}}(\overline{\bm{x}},\bm{\xi})$\\
        $\bm{\tau}^{\mathrm{SGS}}(\overline{\bm{x}},\bm{\xi}) \sim p_{data}(\cdot | \overline{\bm{x}}) $\\
        $\overline{\bm{x}}(0) = \overline{\bm{x}}_0$};
    
    \draw[->,thick] (node cs:name=a) -- (node cs:name=b);
    \draw[->,thick] (node cs:name=b) -- (node cs:name=c);

    \node[draw=none,fill=white] at (2.15,0) {filter};
    \node[draw=none,fill=white] at (6.8,0) {ROM};
    
    \end{tikzpicture}
\end{center}}
\noindent
where, to close the coarse-scale system, we have proposed to model $\bm{\tau}(\overline{\bm{x}},\bm{x}) = \overline{\bm{f}(\bm{x})} - \bm{f}(\overline{\bm{x}})$ as a distribution conditioned on the coarse-scale states $\overline{\bm{x}}$
\begin{align}
    \bm{\tau}(\overline{\bm{x}},\bm{x}) \approx \bm{\tau}^{\mathrm{SGS}} (\overline{\bm{x}},\bm{\xi}) \sim p_{data}(\cdot | \overline{\bm{x}})
\end{align}
The simulation of an SDE can be used for converting a simple distribution $p_{init}$ into a complex distribution $p_{data}$. A diffusion model is described by the SDE
\begin{align}
    \bm{t}_0 &\sim p_{init} = \mathcal{N}(0,I_d) \\
    d\bm{t}_\gamma &= \left[\tilde{u}^\theta_\gamma (\bm{t}_\gamma | \overline{\bm{x}}) + \frac{\sigma_\gamma^2}{2} \tilde{s}^\theta_\gamma (\bm{t}_\gamma | \overline{\bm{x}}) \right] d\gamma + \sigma_\gamma dW_\gamma \\
    \bm{t}_1 &\sim p_{data}(\cdot | \overline{\bm{x}}) = \bm{\tau}^{\mathrm{SGS}}
\end{align}
where, for a fixed choise of guidance scale $w > 1$, we define
\begin{align}
    \tilde{s}^\theta_\gamma (\bm{t}|\overline{\bm{x}}) &= (1-w)s^\theta_\gamma (\bm{t}|\varnothing) + w s^\theta_\gamma (\bm{t}|\overline{\bm{x}}) \\
    \tilde{u}^\theta_\gamma (\bm{t}|\overline{\bm{x}}) &= (1-w)u^\theta_\gamma(\bm{t}|\varnothing) + w u^\theta_\gamma (\bm{t}|\overline{\bm{x}})
\end{align}
The diffusion model consists of a vector field ${u}_\gamma^\theta$, a score network ${s}_\gamma^\theta$, and a fixed diffusion coefficient $\sigma_\gamma$. If we set $\sigma_\gamma = 0$, we recover a flow model. The flow training target $u_\gamma^{\mathrm{target}}$ is the marginal vector field.  To construct a training target $u_\gamma^{\mathrm{target}}$ for a flow model, we choose a conditional probability path $p_\gamma(\bm{t}|\bm{\tau})$ that fulfills $p_0(\cdot|\bm{\tau}) = p_{init}$, $p_1(\cdot|\bm{\tau}) = \delta_{\bm{\tau}}$ (Dirac delta distribution). We consider a Gaussian conditional probability path $p_\gamma (\bm{t}|\bm{\tau}) = \mathcal{N}(\bm{t}; \alpha_\gamma \bm{\tau}, \beta_\gamma^2 I_d)$ where $\alpha_\gamma = \gamma$ and $\beta_\gamma = \sqrt{1-\gamma}$ are monotonic, continuously differentiable functions satisfying $\alpha_1 = \beta_0 = 1$ and $\alpha_0 = \beta_1 = 0$. Guided flow matching consists of training a model $u_\gamma^\theta$ via minimizing the classifier-free guidance (CFG) conditional flow matching (CFM) loss
\begin{align}
    \mathcal{L}^{\mathrm{CFG}}_{\mathrm{CFM}} &= \mathbb{E}_\square ||u^\theta_\gamma (\bm{t}|\overline{\bm{x}}) - u^{\mathrm{target}}_\gamma (\bm{t}|\bm{\tau})||^2 \\
    \square &= (\bm{\tau},\overline{\bm{x}}) \sim p_{data}(\bm{\tau},\overline{\bm{x}}), \gamma \sim \mathrm{Unif}, \bm{t} \sim p_\gamma(\cdot|\bm{\tau}), \text{ replace } \overline{\bm{x}} \text{ with } \varnothing \text{ with probability } \eta
\end{align}
where $u_\gamma^{\mathrm{target}}(\bm{t}|\bm{\tau})$ is the conditional vector field. We train a score network $s_\gamma^\theta$ via CFG denoising score matching (DSM)
\begin{align}
    \mathcal{L}^{\mathrm{CFG}}_{\mathrm{DSM}} &= \mathbb{E}_\square ||s^\theta_\gamma (\bm{t}|\overline{\bm{x}}) - \nabla \log p_\gamma (\bm{t}|\bm{\tau})||^2 \\
    \square &= (\bm{\tau},\overline{\bm{x}}) \sim p_{data}(\bm{\tau},\overline{\bm{x}}), _\gamma \sim \mathrm{Unif}, \bm{t} \sim p_\gamma (\cdot|\bm{\tau}), \text{ replace } \overline{\bm{x}} = \varnothing \text{ with probability } \eta 
\end{align}
where the conditional score distribution is
\begin{align}
    \nabla \log p_\gamma(\bm{t}|\bm{\tau}) &= -\frac{\bm{t} - \alpha_\gamma \bm{\tau}}{\beta_\gamma^2}
\end{align}
See Appendix \ref{sec:cfg_loss} for more details. In practice, $(\bm{\tau},\overline{\bm{x}}) \sim p_{data}(\bm{\tau},\overline{\bm{x}})$ is obtained by sampling a subgrid stress and filtered state from our labeled dataset. The training data includes the exact SGS terms $\bm{\tau}(\overline{\bm{x}},\bm{x}) = \overline{\bm{f}(\bm{x})} - \bm{f}(\overline{\bm{x}})$ and the filtered states $\overline{\bm{x}}$ for $\Delta = 0.01$. The fine-scale trajectory needs to be evolved long enough cover a sufficient amount of the phase-space. At inference time, we can combine $s^\theta_\gamma (\bm{t}|\overline{\bm{x}})$ with $u^\theta_\gamma (\bm{t}|\overline{\bm{x}})$ and sample via
simulating the SDE from $\gamma = 0$ to 1. The goal is for $\bm{t}_1$ to adhere to the guiding variable $\overline{\bm{x}}$. We first use flow matching ($\sigma_\gamma = 0$) with guidance scale $w = 3.0$ and unconditional training probability $\eta = 0.1$. For the flow model, we use an MLP network with 2 hidden layers with 128 neurons per layer and sigmoid linear unit (SiLU) activation functions. Diffusion timestep $d\gamma = 0.00014$ for the diffusion sampling process. The structure is enforced in the sampling process (i.e., linear dynamics result in zero subgrid forcing). Histograms of the initial (top) and learned (bottom) SGS terms for each state are shown in Figure \ref{fig:guided_pdf_flow}. Figure \ref{fig:guided_tau_flow} shows the true and generated SGS terms, while the magnitudes are inaccurate, the transitions are well captured.

\begin{figure}[H]
    \centering
    \includegraphics[width=0.9\textwidth]{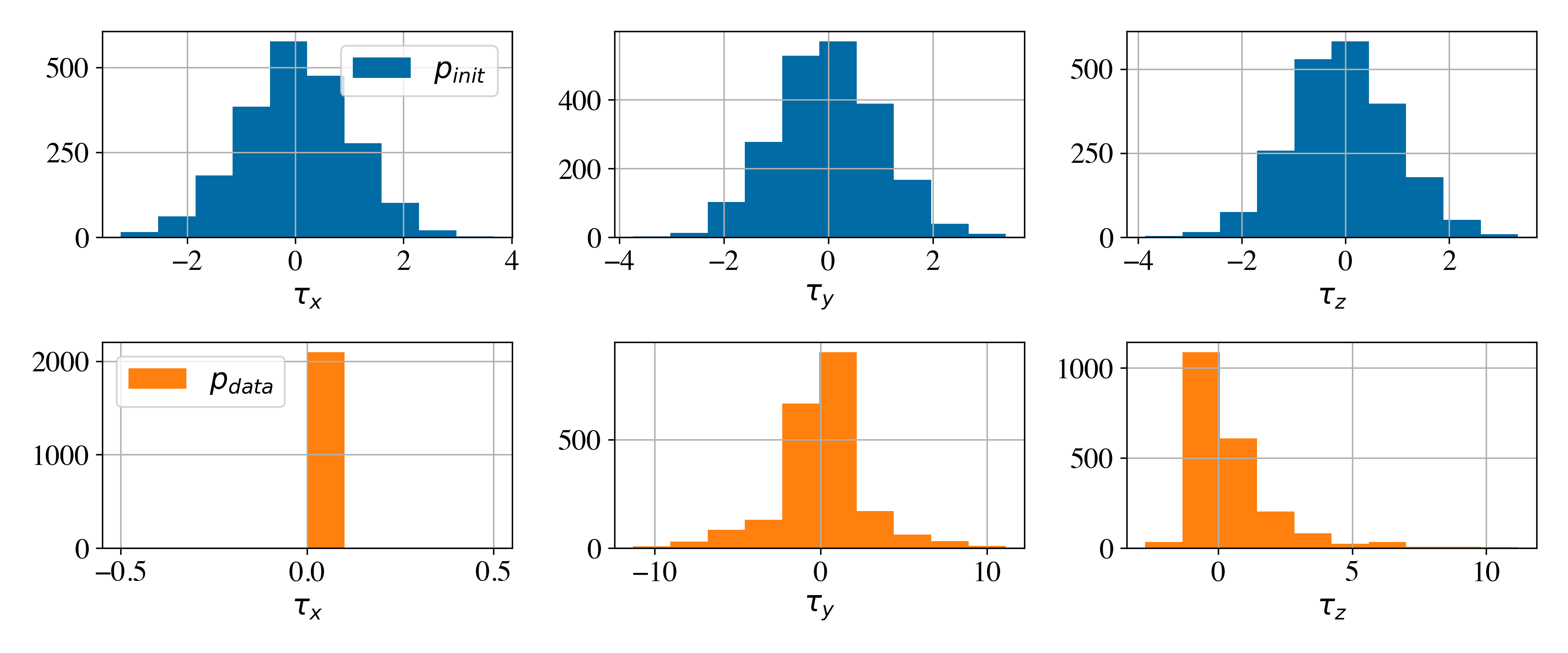}
    \caption{Initial $p_{init}$ and learned $p_{data}$ distributions for guided SGS terms with Langevin dynamics.}
    \label{fig:guided_pdf_sde}
\end{figure}

\begin{figure}[H]
    \centering
    \includegraphics[width=0.9\textwidth]{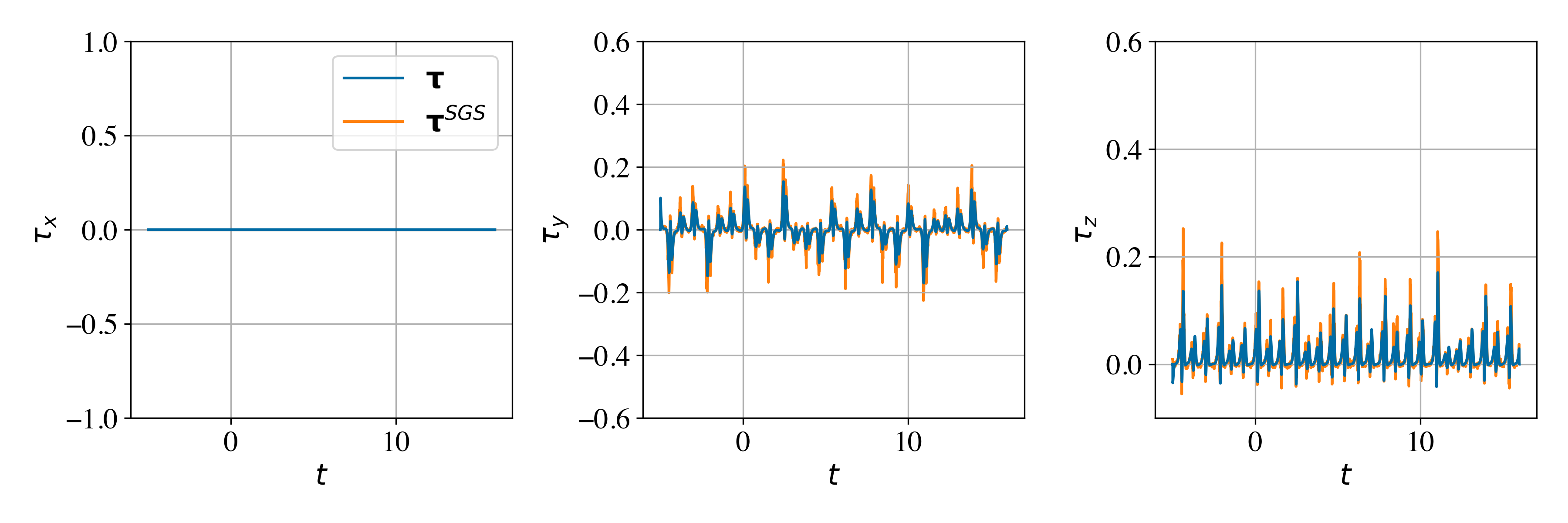}
    \caption{True $\bm{\tau}$ and generated $\bm{\tau}^{\mathrm{SGS}}$ with guided diffusion model.}
    \label{fig:guided_tau_sde}
\end{figure}

\noindent
Now, we add in score matching, with guidance scale $w = 1.5$. For the score model, we use an MLP network with 4 hidden layers with 128 neurons per layer and sigmoid linear unit (SiLU) activation functions. Diffusion constant $\sigma_\gamma = 0.15$ with $d\gamma = 0.0001$ for the diffusion sampling process. Histograms of the initial (top) and learned (bottom) SGS terms for each state are shown in Figure \ref{fig:guided_pdf_sde}. Figure \ref{fig:guided_tau_sde} shows the true and generated SGS terms, where now both the transitions and magnitudes are well captured. These results demonstrate the potential benefits of incorporating stochasticity when closing coarse-scale chaotic systems.

\begin{figure}[H]
    \centering
    \includegraphics[width=0.9\textwidth]{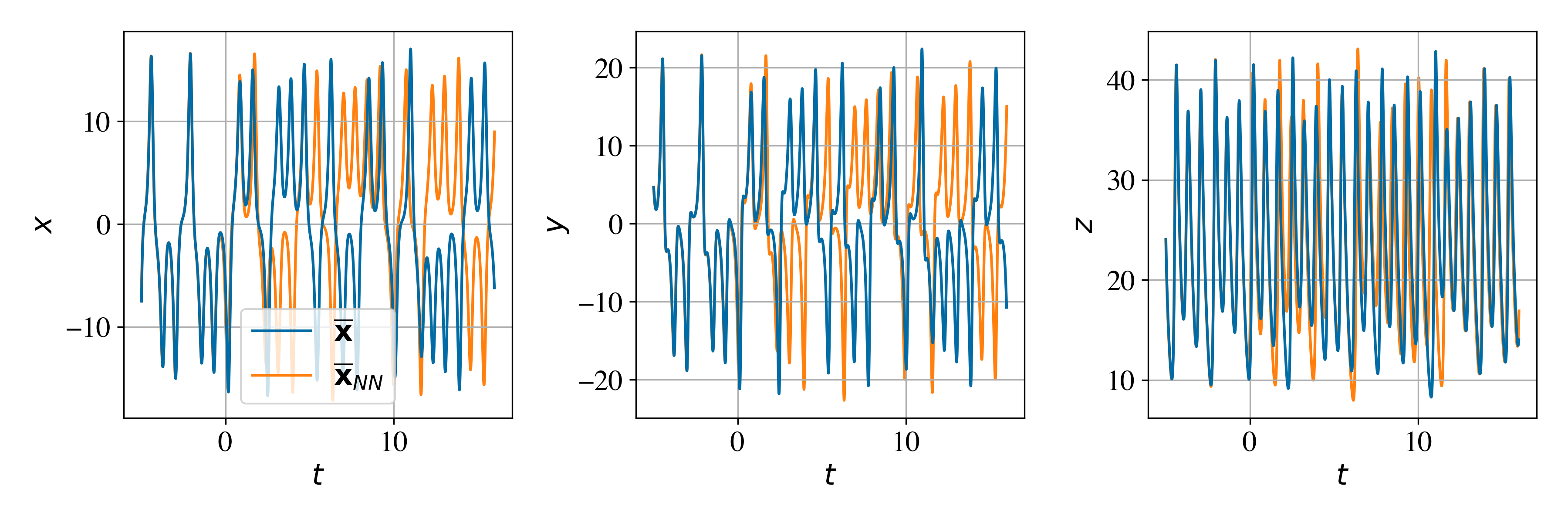}
    \caption{Time evolution of dynamical system with guided SGS model.}
    \label{fig:guided_states_subplots}
\end{figure}

\begin{figure}[H]
    \centering
    \includegraphics[width=0.95\textwidth]{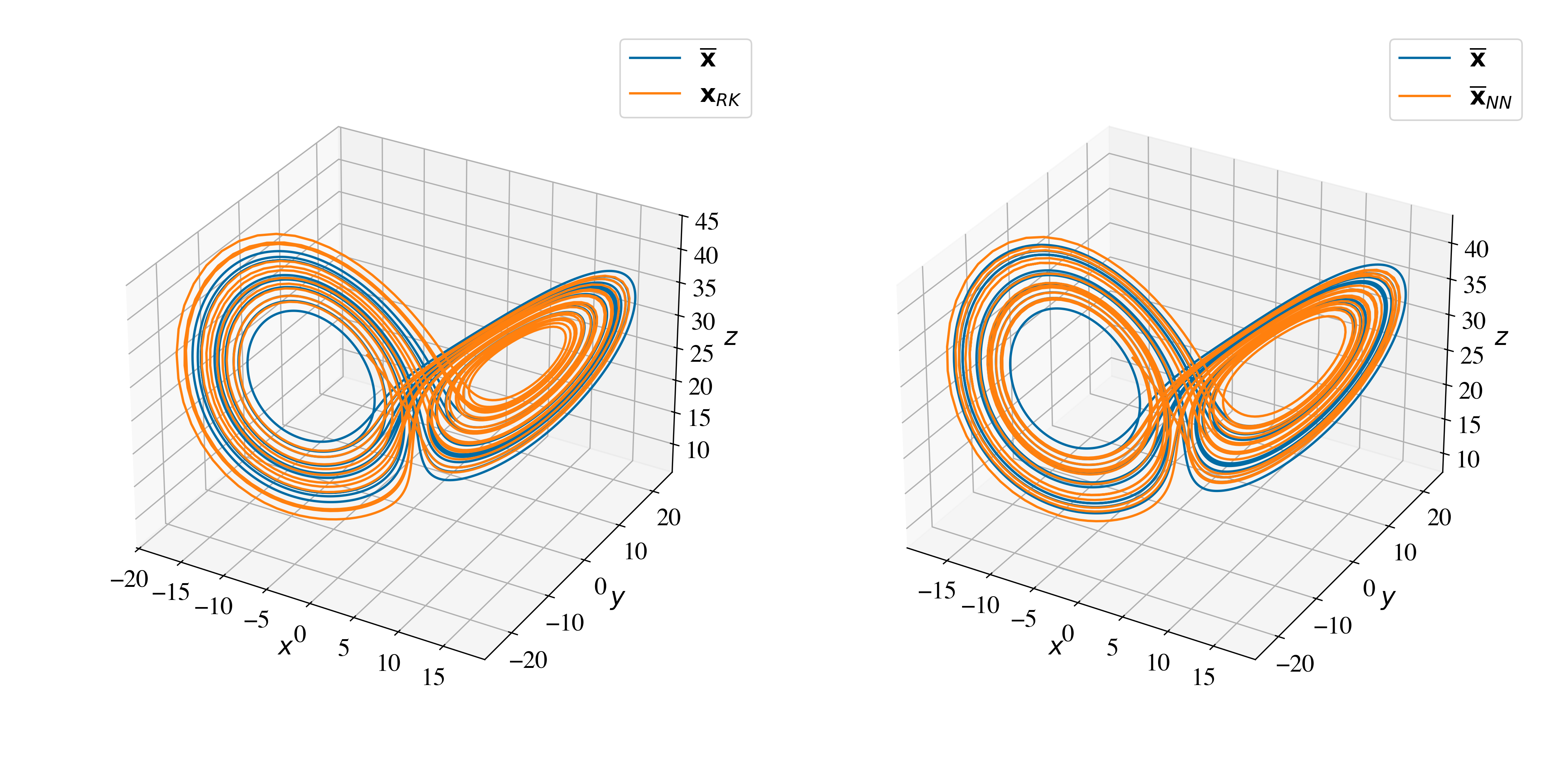}
    \caption{Comparison of (left) deterministic RK $\bm{x}_{RK}$ with timestep $h$ and (right) guided generative method $\overline{\bm{x}}_{NN}$  with timestep $h$ to true filtered states $\overline{\bm{x}}$.}
    \label{fig:guided_states}
\end{figure}

\noindent
Assessing the evolution of the states in Figure \ref{fig:guided_states_subplots} shows realistic transitions between the two focal points of the dynamical system. Figure \ref{fig:guided_states} (right) shows the trajectory using Runge-Kutta (RK) integration with the guided generative SGS model (labeled NN) with timestep $h = \Delta$. Both trajectories are compared to the true $\overline{\bm{x}}$ from filtering the fine-scale trajectory. The RK and NN trajectories start from the same initial condition as the filtered trajectory. The RK trajectory compared with the true $\overline{\bm{x}}$ has Wasserstein distances of $(x,y,z) = (1.81, 1.97, 0.47)$. The trajectory with the generative SGS model compared with the true $\overline{\bm{x}}$ has Wasserstein distances of $(x,y,z) = (0.66, 0.77, 0.10)$. Thus, the trajectory with the generative SGS model is closer than the RK trajectory to the true $\overline{\bm{x}}$ for the same coarse timestep.

\subsection{Linearized stabilization}
Now, consider the following distribution-based stabilization for the linearized dynamics.

\begin{figure}[H]
    \centering
    \includegraphics[width=0.9\textwidth]{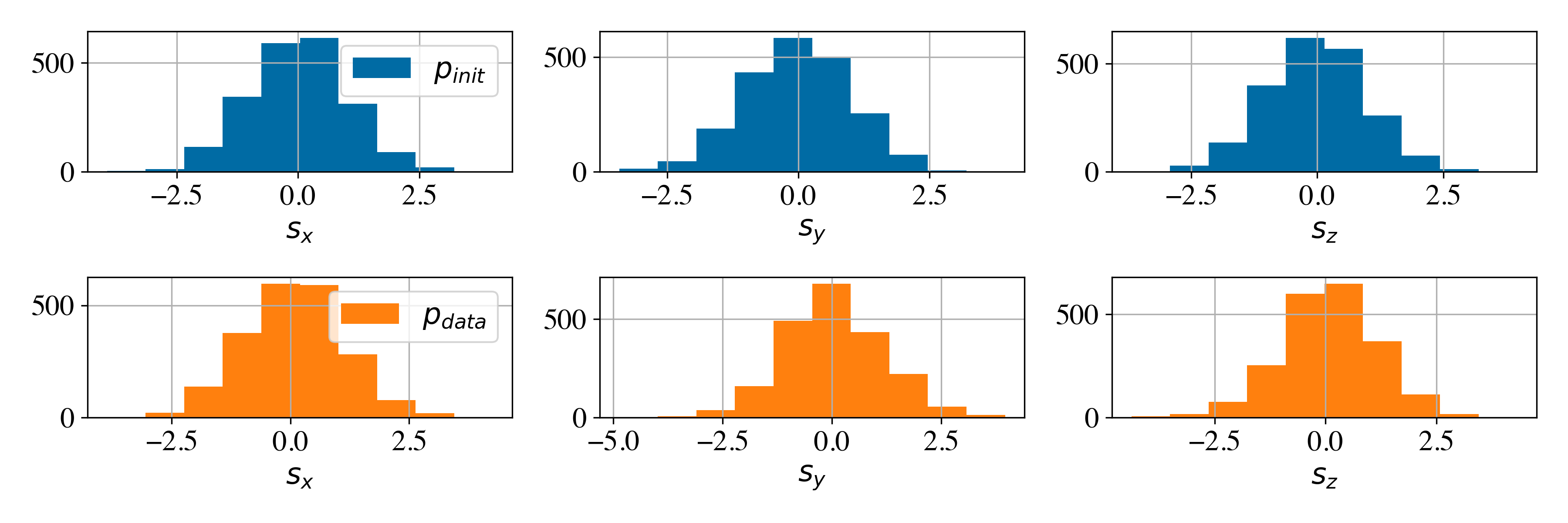}
    \caption{Initial $p_{init}$ and learned $p_{data}$ distributions for guided stabilization terms.}
    \label{fig:guided_pdf_linearized}
\end{figure}

\begin{figure}[H]
    \centering
    \includegraphics[width=0.9\textwidth]{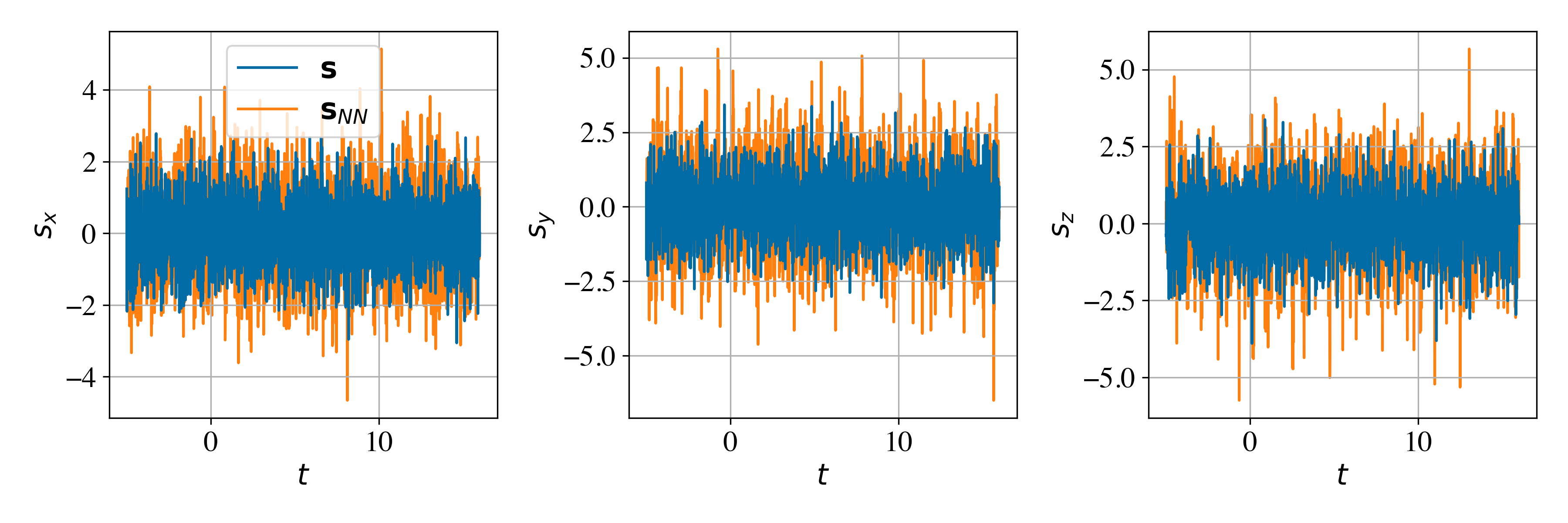}
    \caption{True $\bm{s}$ and generated $\bm{s}_{NN}$ with guided diffusion model.}
    \label{fig:guided_s_linearized}
\end{figure}

{\fontsize{11}{14}
    \linespread{1.75}
    \begin{center}
    \begin{tikzpicture}[
    squarednode/.style={rectangle, rounded corners, draw=teal!25, fill=cyan!5, minimum height=0.5em, minimum width=1.5em},
    ]

    \path (0,0) node(a) [rectangle, rounded corners, draw=teal!25, fill=cyan!5, minimum height=0.5em, minimum width=1.5em, align=left, inner xsep=0.5em, inner ysep=0.5em]
    {   $\bm{x}_t = \bm{f}(\bm{x})$\\
        $\bm{x}(0) = \bm{x}_0$};
    \path (5,0) node(b) [rectangle, rounded corners, draw=teal!25, fill=cyan!5, minimum height=0.5em, minimum width=1.5em, align=left, inner xsep=0.5em, inner ysep=0.5em]
    {   $\tilde{\bm{x}}_t = \bm{f_x}(\bm{x}) \tilde{\bm{x}}$\\
        $\tilde{\bm{x}}(0) = \tilde{\bm{x}}_0$};
    \path (11,0) node(c) [rectangle, rounded corners, draw=teal!25, fill=cyan!5, minimum height=0.5em, minimum width=1.5em, align=left, inner xsep=0.5em, inner ysep=0.5em] 
    {   $\tilde{\bm{x}}_t = \bm{f_x}(\bm{x}) \tilde{\bm{x}} + \bm{s}(\bm{x},\tilde{\bm{x}},\bm{\xi})$\\
    $\bm{s}(\bm{x},\tilde{\bm{x}},\bm{\xi}) \sim p_{data}(\cdot | \bm{x},\tilde{\bm{x}})$\\
        $\tilde{\bm{x}}(0) = \tilde{\bm{x}}_0$};
    
    \draw[->,thick] (node cs:name=a) -- (node cs:name=b);
    \draw[->,thick] (node cs:name=b) -- (node cs:name=c);

    \node[draw=none,fill=white] at (2.4,0) {linearize};
    \node[draw=none,fill=white] at (7.5,0) {stabilize};
    
    \end{tikzpicture}
\end{center}}
\noindent
where, to stabilize the linearized system, we propose to model $\bm{s}(\bm{x},\tilde{\bm{x}},\bm{\xi})$ as a distribution conditioned on the nonlinear $\bm{x}$ and linearized $\tilde{\bm{x}}$ states
\begin{align}
    \bm{s}(\bm{x},\tilde{\bm{x}},\bm{\xi}) \sim p_{data}(\cdot | \bm{x},\tilde{\bm{x}})
\end{align}
The training data is collected following Algorithm \ref{alg:train2} to generate the following distribution $p_{data}$ for $\bm{s}$
\begin{align}
    \bm{s}_n^k = \tilde{\bm{x}}_{n+1}^k - \frac{1}{\Delta t}( \tilde{\bm{x}}_n^k + \bm{f}_{\bm{x}}(\bm{x}_n) \tilde{\bm{x}}_n^k \Delta t) \sim p_{data}
\end{align}
Flow (with $\eta = 0.2$) and score (with $\eta = 0.5$) models are trained to minimize CFG matching training objectives. At inference time, we can combine score $\tilde{s}^\theta_\gamma (\bm{t}|{\bm{x}},\tilde{\bm{x}})$ with flow $\tilde{u}^\theta_\gamma (\bm{t}|{\bm{x}},\tilde{\bm{x}})$ models and sample via
\begin{align}
    \bm{t}_0 &\sim p_{init} = \mathcal{N}(0,I_d) \\
    d\bm{t}_\gamma &= \left[\tilde{u}^\theta_\gamma (\bm{t}_\gamma | {\bm{x}},\tilde{\bm{x}}) + \frac{\sigma_\gamma^2}{2} \tilde{s}^\theta_\gamma (\bm{t}_\gamma | {\bm{x}},\tilde{\bm{x}}) \right] d\gamma + \sigma_\gamma dW_\gamma \\
    \bm{t}_1 &\sim p_{data}(\cdot | {\bm{x}},\tilde{\bm{x}}) = \bm{s}
\end{align}
simulating the SDE from $\gamma = 0$ to 1. The goal is for $\bm{t}_1$ to adhere to the guiding variables ${\bm{x}}$ and $\tilde{\bm{x}}$. For the flow and score models, we use MLP networks with 4 hidden layers with 128 neurons per layer each and sigmoid linear unit (SiLU) activation functions. Diffusion constant $\sigma_\gamma  = 0.5$ with $d\gamma = 0.00014$ for the diffusion sampling process with guidance scale $w = 0.1$. Histograms of the initial (top) and learned (bottom) stabilization terms for each state are shown in Figure \ref{fig:guided_pdf_linearized}. Figure \ref{fig:guided_s_linearized} shows the true and generated $\bm{s}$ terms.

\begin{figure}[H]
    \centering
    \includegraphics[width=0.9\textwidth]{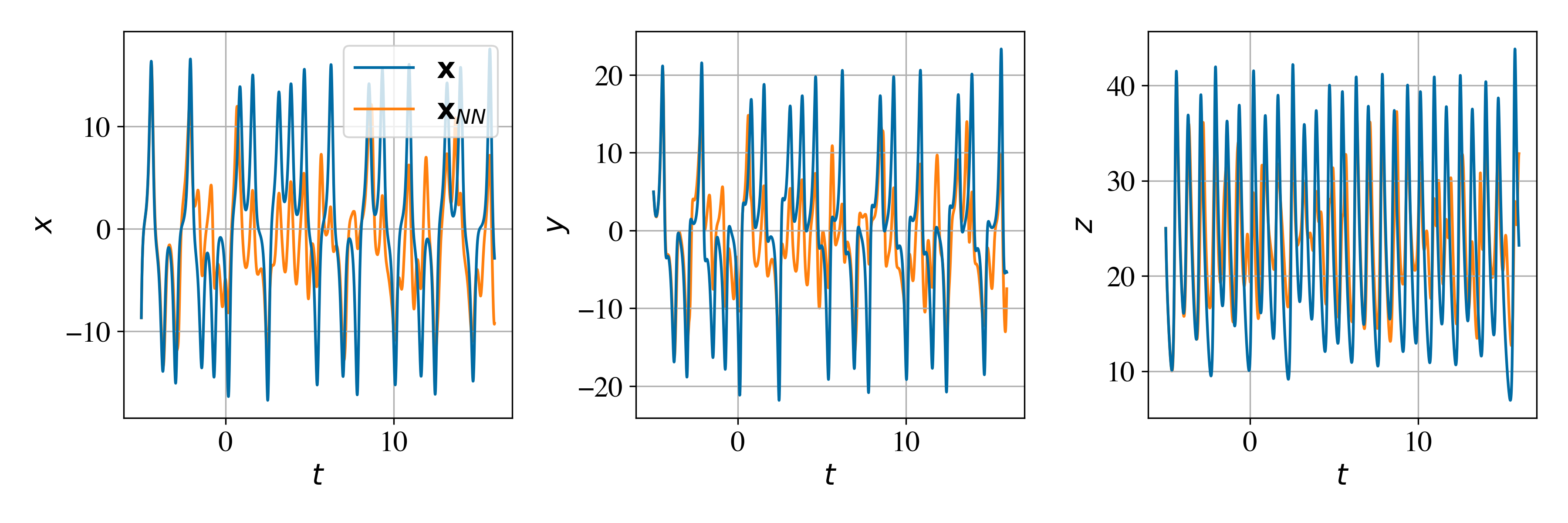}
    \caption{Time evolution of dynamical system with guided stabilization model.}
    \label{fig:guided_states_sde_linearized_subplots}
\end{figure}

\noindent
Assessing the evolution of the states $\bm{x}_{NN} = \bm{x} + \tilde{\bm{x}}_{NN}$ in Figure \ref{fig:guided_states_sde_linearized_subplots} shows stabilizing transitions between the two focal points of the dynamical system, without experiencing exponential blow-up present in the linearized system. This demonstrates the applicability of using a probabilistic approach for stabilizing an unstable chaotic system.  

\section{Quadratic diffusion model}

\begin{figure}[H]
    \centering
    \includegraphics[width=0.9\textwidth]{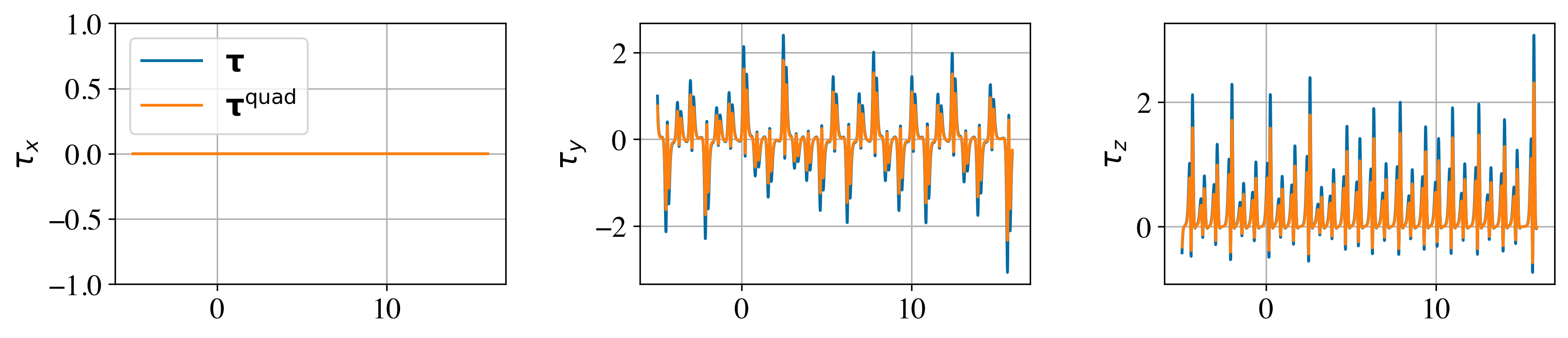}
    \caption{Analytical quadratic model for subgrid states using true fluctuations.}
    \label{fig:exact_sgs_quad}
\end{figure}

We can derive an analytical quadratic model for the subgrid stresses. See Appendix \ref{sec:quad_model_derivation} for the full derivation. For some filter width $\Delta$, the subgrid term can be modeled by
\begin{align}
    &\bm{\tau}^{\mathrm{quad}}(\overline{\bm{x}},\bm{x}) = \frac{1}{2} \bm{x}'^T \left[ \bm{H} + \frac{\Delta^2}{12} \bm{J}^T \bm{H} \bm{J} \right] \bm{x}' \\
    &\bm{x}' = \bm{x} - \overline{\bm{x}}
\end{align}
where the analytical expression requires access to the true fluctuations, and thus, the fine-scale states. $\bm{J}$ and $\bm{H}$ are the Jacobian and Hessian matrices evaluated at the filtered states. For a filter width of $\Delta = 0.04$, the true and modeled SGS terms using the true fluctuations in the analytical quadratic model are shown in Figure \ref{fig:exact_sgs_quad}. These results suggest that capturing the contribution of the fluctuations to the dynamics can yield an accurate subgrid model in the form of the analytical quadratic expression.

\begin{figure}[H]
    \centering
    \includegraphics[width=0.9\textwidth]{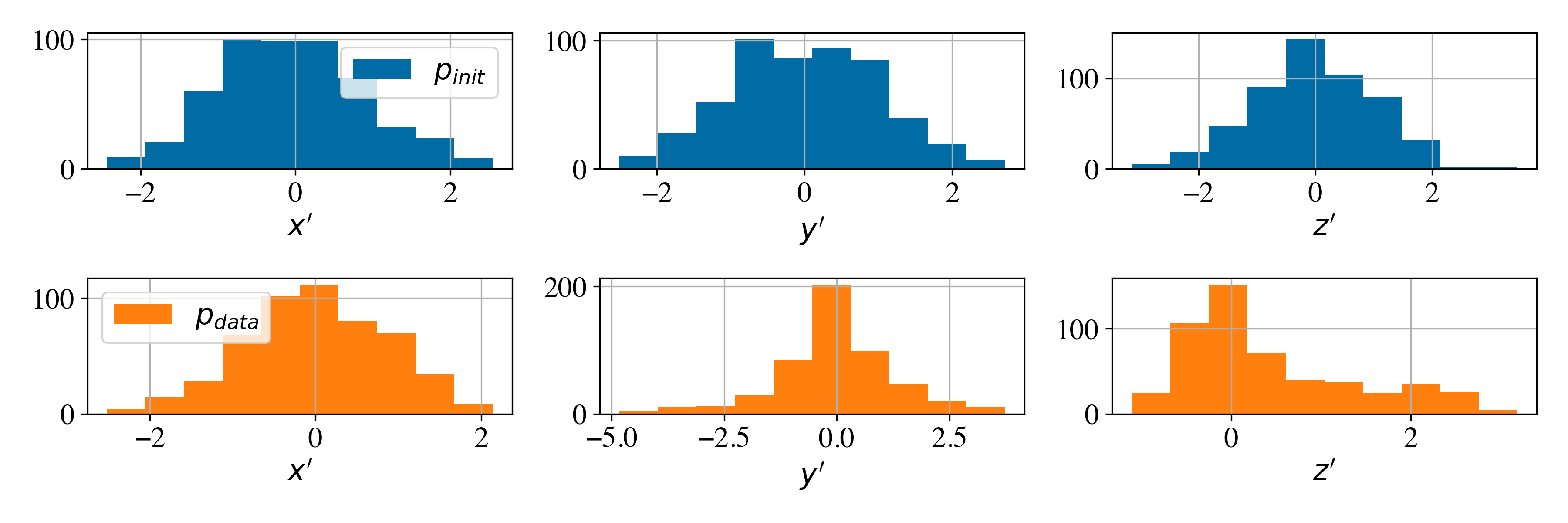}
    \caption{Initial $p_{init}$ and learned $p_{data}$ distributions for guided fluctuations.}
    \label{fig:guided_pdf_xp}
\end{figure}

\begin{figure}[H]
    \centering
    \includegraphics[width=0.9\textwidth]{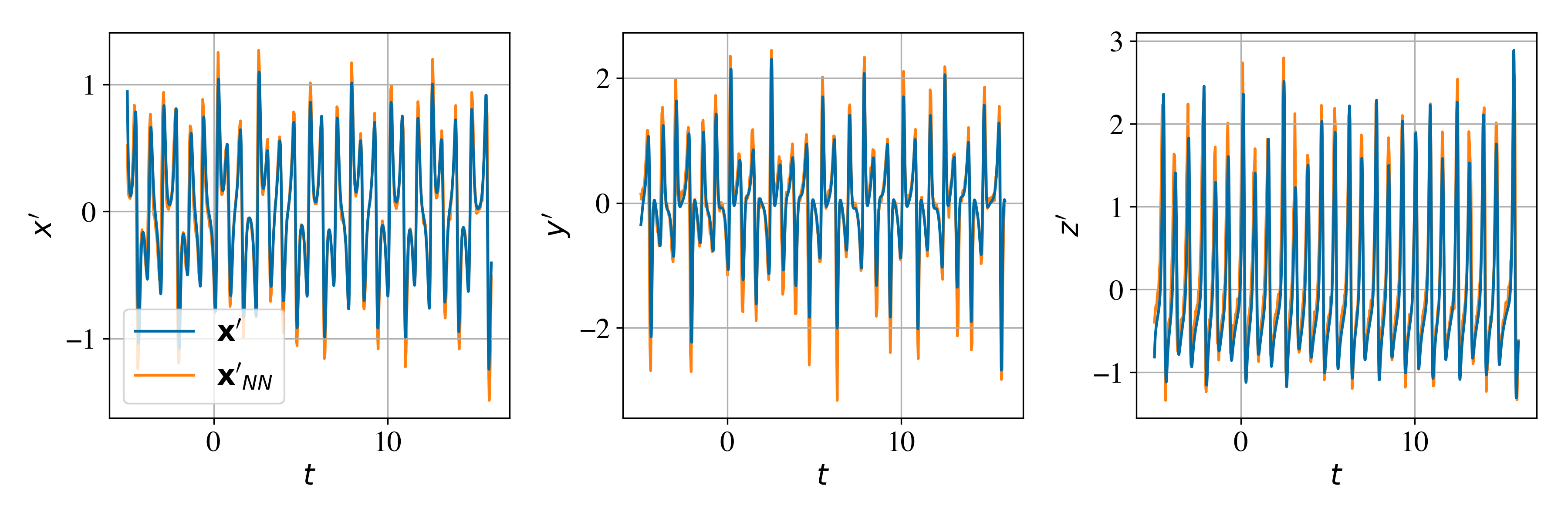}
    \caption{True $\bm{x}'$ and generated $\bm{x}'_{\mathrm{NN}}$ with guided diffusion model.}
    \label{fig:guided_xp}
\end{figure}

\begin{figure}[H]
    \centering
    \includegraphics[width=0.9\textwidth]{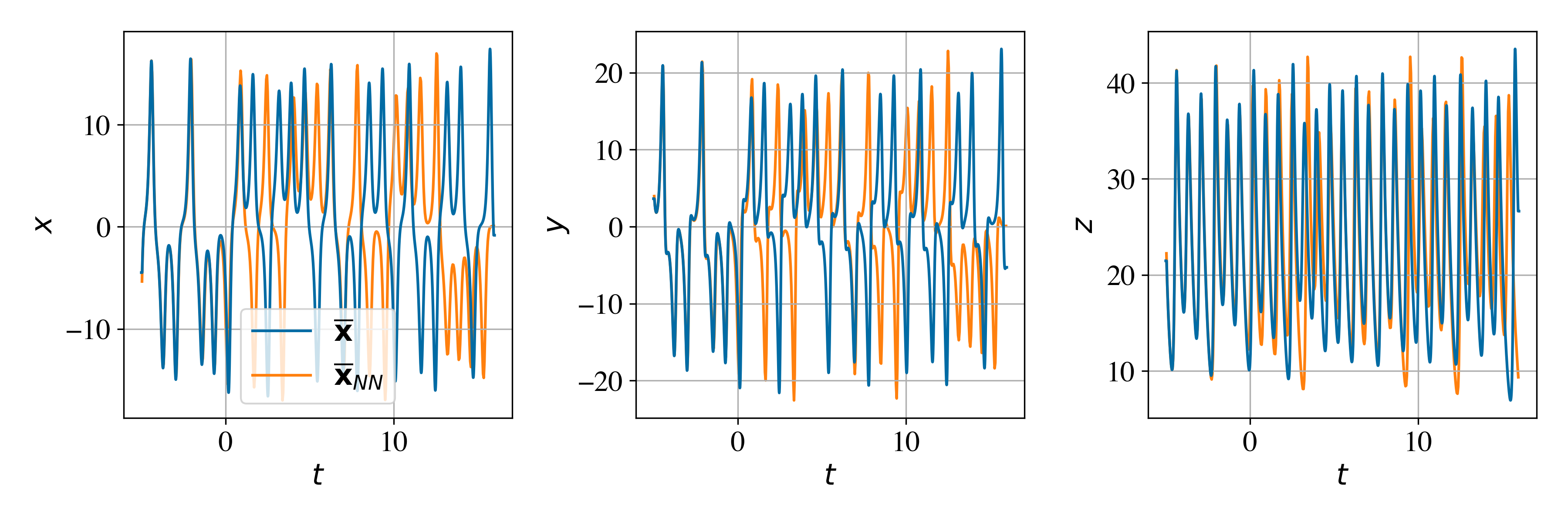}
    \caption{Time evolution of dynamical system with guided quadratic diffusion SGS model.}
    \label{fig:guided_states_sde_quad_subplots}
\end{figure}

\begin{figure}[H]
    \centering
    \includegraphics[width=0.95\textwidth]{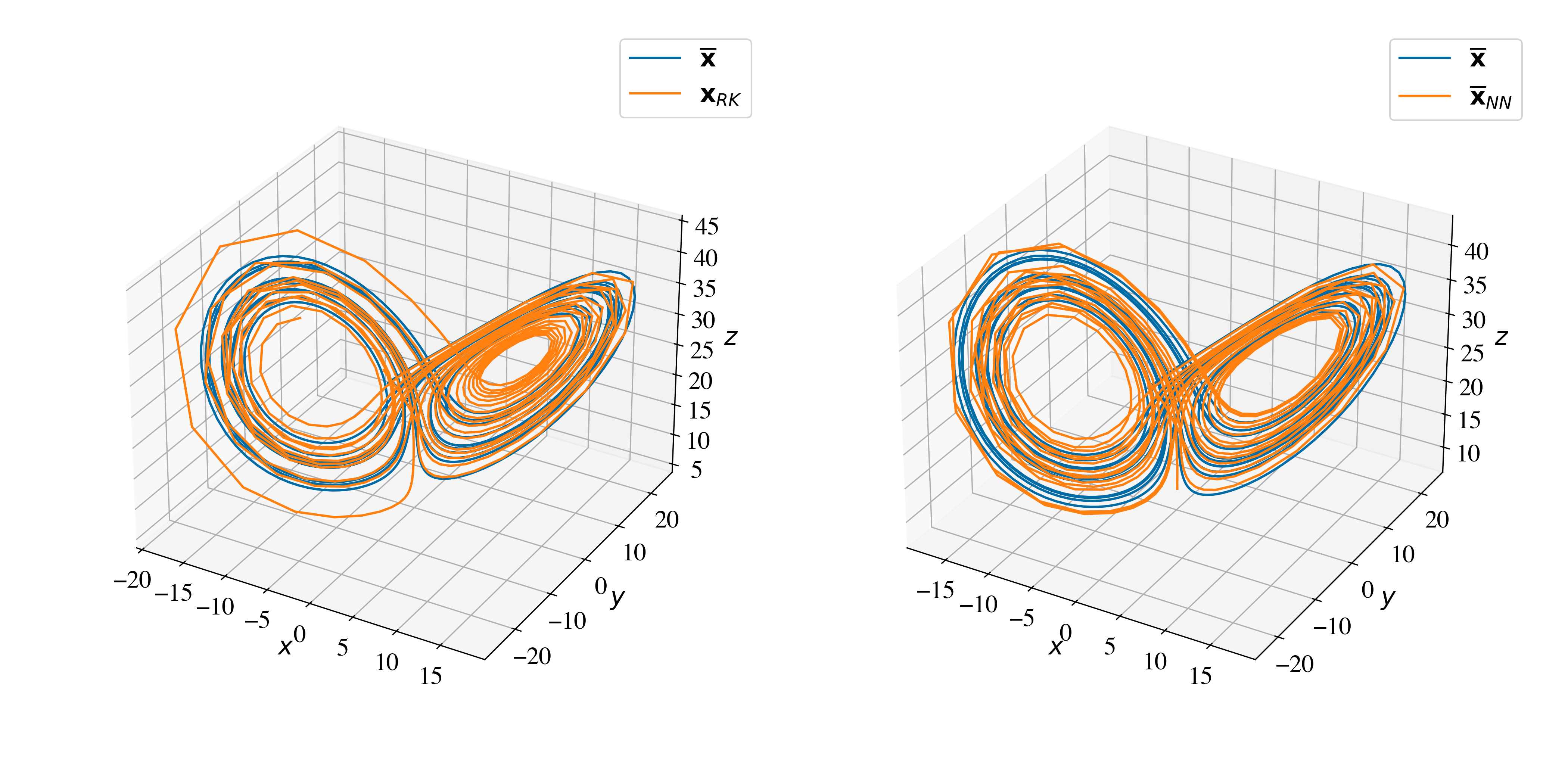}
    \caption{Comparison of (left) deterministic Runge-Kutta $\bm{x}_{RK}$ with timestep $h$ and (right) guided quadratic diffusion SGS model $\overline{\bm{x}}_{NN}$ with timestep $h$ to true filtered states $\overline{\bm{x}}$.}
    \label{fig:guided_states_quad}
\end{figure}

\noindent
We propose using a guided diffusion model for the fluctuations for the quadratic model
{\fontsize{11}{14}
    \linespread{1.75}
    \begin{center}
    \begin{tikzpicture}[
    squarednode/.style={rectangle, rounded corners, draw=teal!25, fill=cyan!5, minimum height=0.5em, minimum width=1.5em},
    ]

    \path (0,0) node(a) [rectangle, rounded corners, draw=teal!25, fill=cyan!5, minimum height=0.5em, minimum width=1.5em, align=left, inner xsep=0.5em, inner ysep=0.5em]
    {   $\bm{x}_t = \bm{f}(\bm{x})$\\
        $\bm{x}(0) = \bm{x}_0$};
    \path (4,0) node(b) [rectangle, rounded corners, draw=teal!25, fill=cyan!5, minimum height=0.5em, minimum width=1.5em, align=left, inner xsep=0.5em, inner ysep=0.5em]
    {   $\overline{\bm{x}}_t = \overline{\bm{f}(\bm{x})}$\\
        $\overline{\bm{x}}(0) = \overline{\bm{x}}_0$};
    \path (10.5,0) node(c) [rectangle, rounded corners, draw=teal!25, fill=cyan!5, minimum height=0.5em, minimum width=1.5em, align=left, inner xsep=0.5em, inner ysep=0.5em] 
    {   $\overline{\bm{x}}_t = \bm{f}(\overline{\bm{x}}) + \bm{\tau}^{\mathrm{quad}}(\overline{\bm{x}},\bm{x})$\\
        $\bm{\tau}^{\mathrm{quad}}(\overline{\bm{x}},\bm{x}) = \frac{1}{2} \bm{x}'^T \left[ \bm{H} + \frac{\Delta^2}{12} \bm{J}^T \bm{H} \bm{J} \right] \bm{x}'$\\
        $\bm{x}'(\overline{\bm{x}},\bm{\xi}) = \bm{x} - \overline{\bm{x}} \sim p_{data}(\cdot | \overline{\bm{x}})$\\
        $\overline{\bm{x}}(0) = \overline{\bm{x}}_0$};
    
    \draw[->,thick] (node cs:name=a) -- (node cs:name=b);
    \draw[->,thick] (node cs:name=b) -- (node cs:name=c);

    \node[draw=none,fill=white] at (2,0) {filter};
    \node[draw=none,fill=white] at (6,0) {ROM};
    
    \end{tikzpicture}
\end{center}}
\noindent
where we propose to sample the fluctuations from a distribution conditioned on the filtered states
\begin{align}
    \bm{x}'(\overline{\bm{x}},\bm{\xi}) = \bm{x} - \overline{\bm{x}} \sim p_{data}(\cdot | \overline{\bm{x}})
\end{align}
The fine-scale trajectory is simulated with $\Delta t = 0.01$ using a second-order discontinuous Galerkin integration scheme. A filter width of $\Delta = 0.04$ is applied to acquire the coarse-scale states. For the flow model, we use an MLP network with 2 hidden layers with 128 neurons per layer and SiLU activation function. For the score model, we use an MLP network with 4 hidden layers with 128 neurons per layer and SiLU activation function. For CFG, we set probability $\eta = 0.1$ and guidance scale $w = 1.5$. Diffusion constant $\sigma_\gamma  = 0.1$ with $d\gamma = 0.0001$ for the diffusion sampling process. Histograms of the initial (top) and learned (bottom) fluctuations for each state are shown in Figure \ref{fig:guided_pdf_xp}. Figure \ref{fig:guided_xp} shows the true and generated (NN) fluctuations. Assessing the evolution of the states in Figure \ref{fig:guided_states_sde_quad_subplots} shows more realistic transitions between the two focal points of the dynamical system, without stalling at any point along the trajectory. Figure \ref{fig:guided_states_quad} (right) shows the trajectory using RK integration with the quadratic diffusion SGS model (NN) with timestep $h = \Delta$.

\begin{figure}[H]
    \centering
    \includegraphics[width=0.9\textwidth]{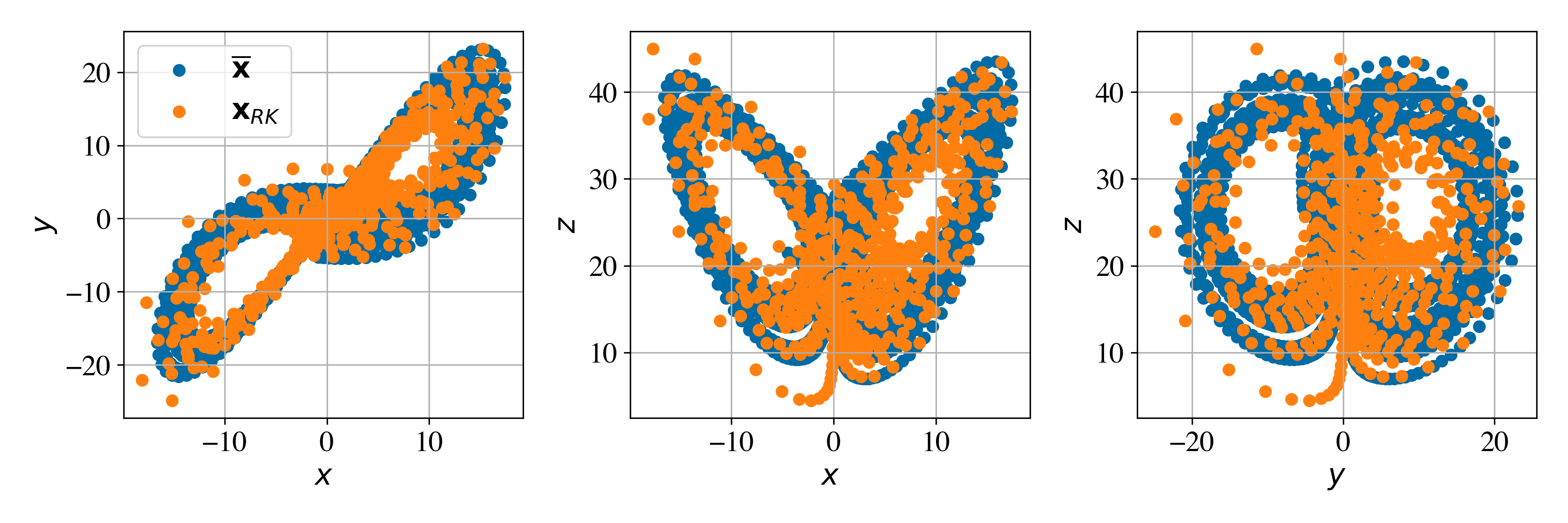}
    \includegraphics[width=0.9\textwidth]{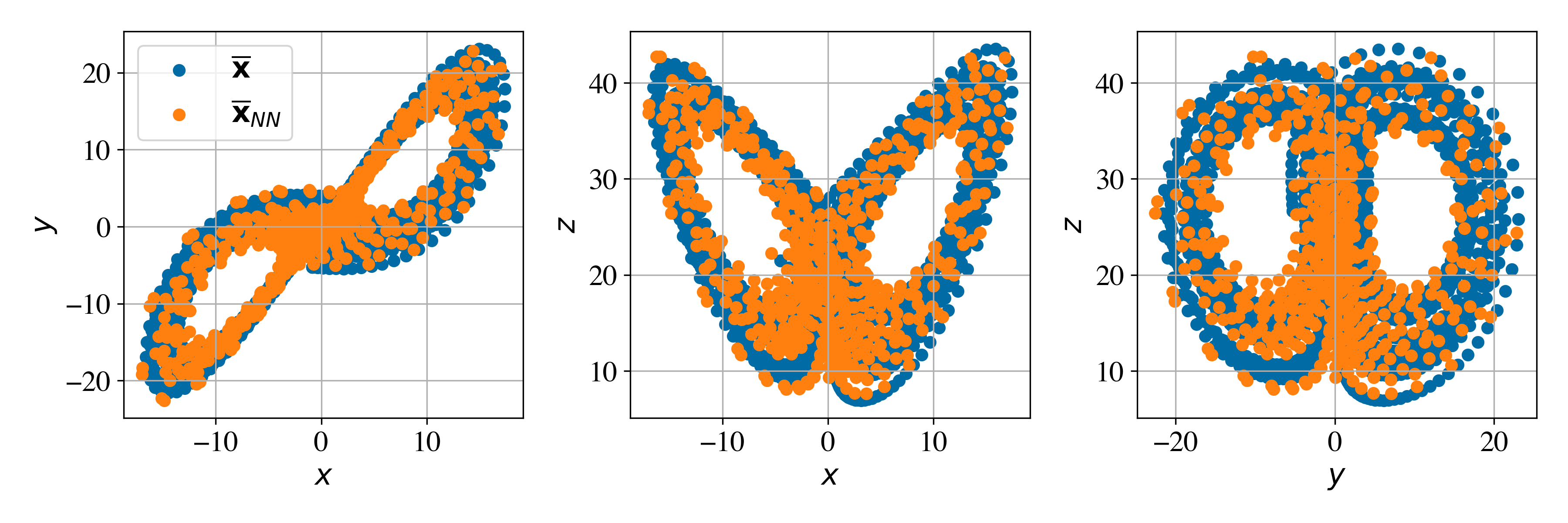}
    \caption{(Top) RK and (bottom) NN projected trajectories compared to true coarse-scale attractor.}
    \label{fig:projected}
\end{figure}

\noindent
Both trajectories are compared to the true $\overline{\bm{x}}$ from filtering the fine-scale trajectory. The RK and NN trajectories start from the same initial condition as the filtered trajectory. The trajectory with the generative SGS model is closer than the RK trajectory to the true $\overline{\bm{x}}$ for the same coarse timestep. The quadratic diffusion model permits use of a far coarser timestep, $h = 4\Delta t$ in this case, to recover the resolved states. These results suggest potential benefits in incorporating theoretical models with generative approaches for closing chaotic systems. Quantitative comparison of two trajectories can be measured using the Hellinger distance, $D_{Hell}$, calculated as the difference between probability densities
\begin{align}
    (D_{Hell})^2 = \frac{1}{2} \int \left( \sqrt{p(\bm{x})} - \sqrt{q(\bm{x})} \right)^2 \,d\bm{x}
\end{align}
where $p(\bm{x})$ and $q(\bm{x})$ are simulated and true probability density functions (pdfs). The smaller the measure, the closer the simulated pdf is to the true pdf \cite{christensen_simulating_2015}. Table \ref{tab:hell} shows the calculated Hellinger distances for projected trajectories, shown in Figure \ref{fig:projected}. 

\begin{table}[H]
    \centering
    \begin{tabular}{|c|c|c|c|c|}\hline
         & $x-y$ & $x-z$ & $y-z$ \\ \hline\hline
    RK & 0.35 & 0.41 & 0.43 \\ \hline
    NN & 0.27 & 0.38 & 0.42 \\ \hline
    \end{tabular}
    \caption{Hellinger distances $D_{Hell}$ for RK and NN trajectories compared to true coarsened filtered attractor.}
    \label{tab:hell}
\end{table}

\begin{figure}[H]
    \centering
    \includegraphics[width=0.9\textwidth]{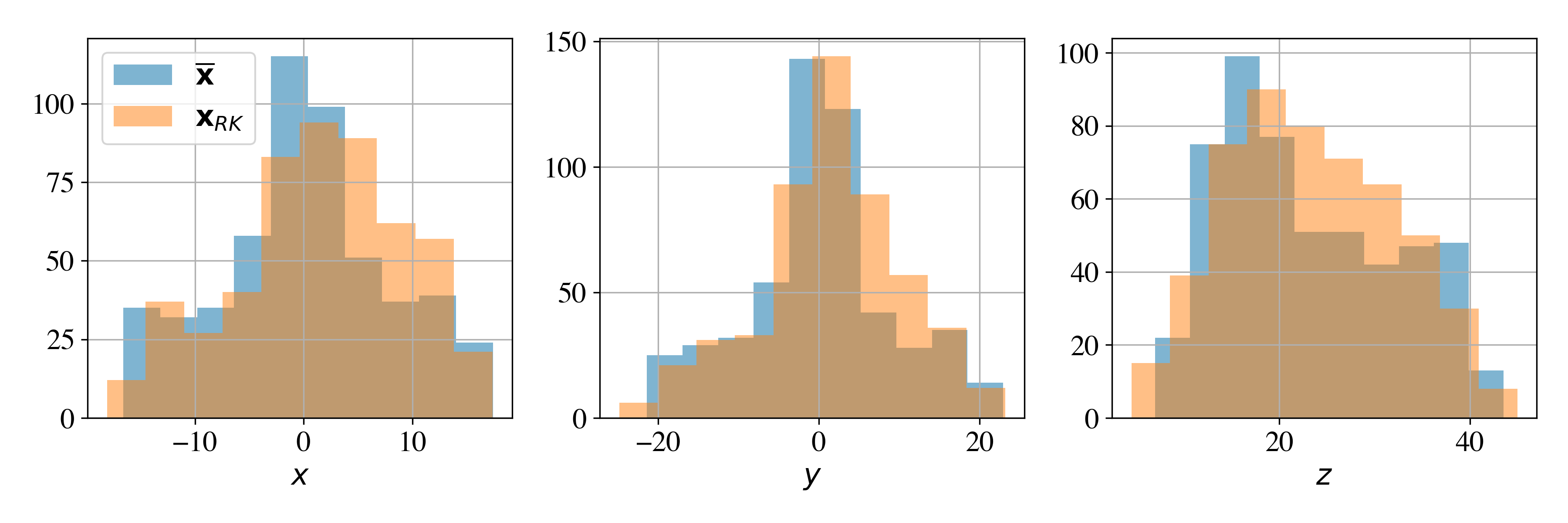}
    \includegraphics[width=0.9\textwidth]{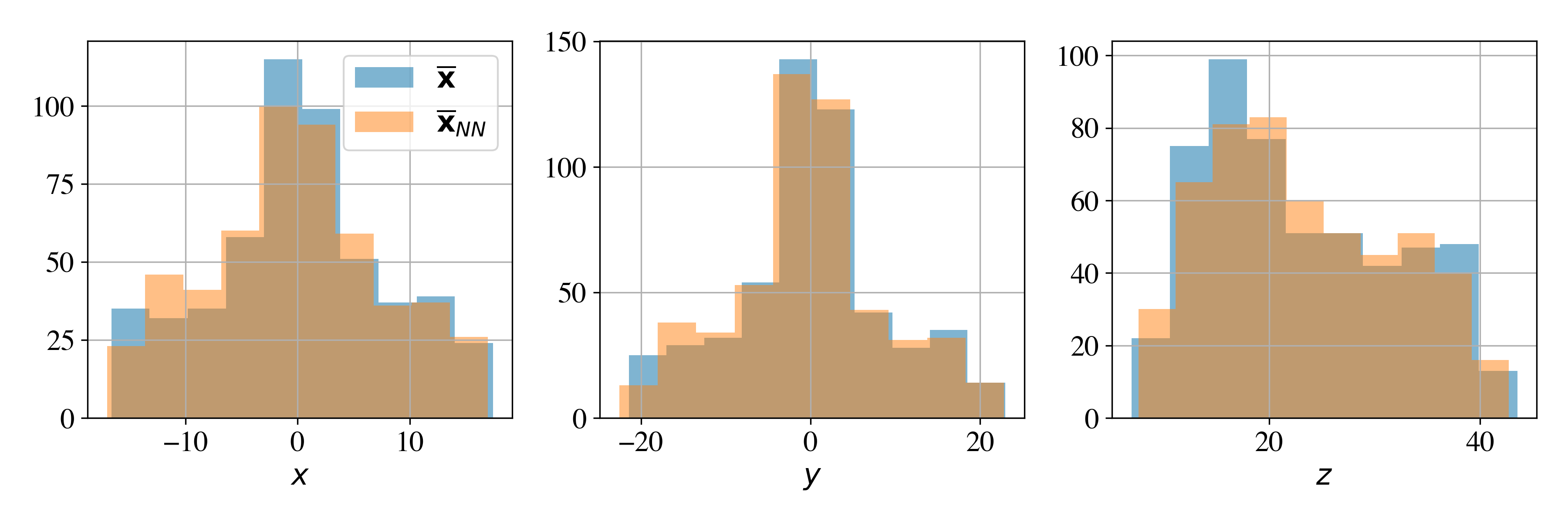}
    \caption{(Top) RK and (bottom) NN distributions compared to true coarsened filtered distributions.}
    \label{fig:hist}
\end{figure}

\noindent
Lastly, comparison of the histograms of each state in Figure \ref{fig:hist} for the RK and NN trajectories compared with the true coarsened filtered trajectory demonstrates better agreement with the guided quadratic diffusion SGS model, quantitatively measured by the Wasserstein distance in Table \ref{tab:wass}.

\begin{table}[H]
    \centering
    \begin{tabular}{|c|c|c|c|c|}\hline
         & $x$ & $y$ & $z$ \\ \hline\hline
    RK & 1.43 & 1.45 & 1.03 \\ \hline
    NN & 0.49 & 0.50 & 0.49 \\ \hline
    \end{tabular}
    \caption{Wasserstein distances $W_1$ for RK and NN distributions compared to true coarsened filtered distribution.}
    \label{tab:wass}
\end{table}

\section{Conclusions}

The results presented demonstrate the applicability of generative SDEs for accurate and stable closure modeling of chaotic systems. We demonstrate that, with the addition of a stochastic generative SGS model, we are able to coarsen the timestep without having significant effect on the accuracy. We also show increased accuracy when combining generative approaches with derived models based on theory, rather than directly generating unresolved dynamics. These results demonstrate the stability and accuracy benefits when incorporating stochastic and generative approaches for closure modeling. Future work includes quantifying rigorous stability guarantees, expanding on existing studies in stochastic stability for linear systems.

Further, we aim to extend the formulation toward closure modeling for chaotic partial differential equations (PDEs). The formulation is conceptually similar to accounting for memory in a dynamical system, as in Mori-Zwanzig (MZ) temporal coarsening \cite{mori_transport_1965, zwanzig_memory_1961, chorin_optimal_2002}. Temporal filtering results in timescales in the ground truth that are not resolved. These scales can be viewed as unresolved variables that have memory-like effects on the resolved variables. In this work, we consider this as adjusting the fluctuations instantaneously to the current filtered state. In particular, the diffusion model conditioned on the filtered states gives the contribution from the fluctuations on the dynamics. More elaborate memory models, such as conditioning on previous filtered states (i.e., longer memory), will be considered in future work.

\section*{Acknowledgments}

The authors acknowledge Renato Trono Figueras for his contributions to initial model development. E. Williams acknowledges Peter Holderrieth and Ezra Erives for their joint instruction of MIT 6.S184: Generative AI with Stochastic Differential Equations \cite{flowsanddiffusions2025}.

This material is based upon work supported by the U.S. Department of Energy, Office of Science, Office of Advanced Scientific Computing Research, Department of Energy Computational Science Graduate Fellowship under Award Number DE-SC0023112. This report was prepared as an account of work sponsored by an agency of the United States Government. Neither the United States Government nor any agency thereof, nor any of their employees, makes any warranty, express or implied, or assumes any legal liability or responsibility for the accuracy, completeness, or usefulness of any information, apparatus, product, or process disclosed, or represents that its use would not infringe privately owned rights. Reference herein to any specific commercial product, process, or service by trade name, trademark, manufacturer, or otherwise does not necessarily constitute or imply its endorsement, recommendation, or favoring by the United States Government or any agency thereof. The views and opinions of authors expressed herein do not necessarily state or reflect those of the United States Government or any agency thereof.

\appendix

\section{Classifier-free guidance training objectives}
\label{sec:cfg_loss}

We have the guided conditional flow matching objective
\begin{align}
    \mathcal{L}^{\mathrm{guided}}_{\mathrm{CFM}} &= \mathbb{E}_\square ||u^\theta_\gamma (\bm{t}|\overline{\bm{x}}) - u^{\mathrm{target}}_\gamma (\bm{t}|\bm{\tau})||^2 \\
    \square &= (\bm{\tau},\overline{\bm{x}}) \sim p_{data}(\bm{\tau},\overline{\bm{x}}), \gamma \sim \mathrm{Unif}, \bm{t} \sim p_\gamma(\cdot|\bm{\tau})
\end{align}
This objective must be amended to account for the possibility of the ``label'' not being in the training set (i.e., $\overline{\bm{x}} = \varnothing)$ \cite{ho2022classifierfreediffusionguidance}. For some guidance scale $w > 1$, define
\begin{align}
    \tilde{u}_\gamma (\bm{t}|\overline{\bm{x}}) = (1 - w) u^{\mathrm{target}}_\gamma (\bm{t}|\varnothing) + w u^{\mathrm{target}}_\gamma (\bm{t}|\overline{\bm{x}})
\end{align}
We can treat $u^{\mathrm{target}}_\gamma (\bm{t})$ as $u^{\mathrm{target}}_\gamma (\bm{t}|\overline{\bm{x}})$, where $\overline{\bm{x}} = \varnothing$ denotes the absence of conditioning. The label set is augmented with a new, additional $\varnothing$ label, so that $\overline{\bm{x}} \in \mathcal{Y}' \triangleq \{\mathcal{Y},\varnothing\}$. Then we have the classifier-free guidance (CFG) conditional flow matching training objective
\begin{align}
    \mathcal{L}^{\mathrm{CFG}}_{\mathrm{CFM}} &= \mathbb{E}_\square ||u^\theta_\gamma (\bm{t}|\overline{\bm{x}}) - u^{\mathrm{target}}_\gamma (\bm{t}|\bm{\tau})||^2 \\
    \square &= (\bm{\tau},\overline{\bm{x}}) \sim p_{data}(\bm{\tau},\overline{\bm{x}}), \gamma \sim \mathrm{Unif}, \bm{t} \sim p_\gamma(\cdot|\bm{\tau}), \text{ replace } \overline{\bm{x}} \text{ with } \varnothing \text{ with probability } \eta
\end{align}
Analogously, for the score network, define the classifier-free guided score $\tilde{s}_\gamma (\bm{t}|\overline{\bm{x}})$ by
\begin{align}
    \tilde{s}_\gamma (\bm{t}|\overline{\bm{x}}) = (1-w) \nabla \log p_\gamma (\bm{t} | \varnothing) + w \nabla \log p_\gamma (\bm{t}|\overline{\bm{x}})
\end{align}
And we arrive at the CFG DSM objective
\begin{align}
    \mathcal{L}^{\mathrm{CFG}}_{\mathrm{DSM}} &= \mathbb{E}_\square ||s^\theta_\gamma (\bm{t}|\overline{\bm{x}}) - \nabla \log p_\gamma (\bm{t}|\bm{\tau})||^2 \\
    \square &= (\bm{\tau},\overline{\bm{x}}) \sim p_{data}(\bm{\tau},\overline{\bm{x}}), \gamma \sim \mathrm{Unif}, \bm{t} \sim p_\gamma(\cdot|\bm{\tau}), \text{ replace } \overline{\bm{x}} = \varnothing \text{ with probability } \eta 
\end{align}

\section{Analytical quadratic model derivation}
\label{sec:quad_model_derivation}

We can derive an analytical quadratic model for the subgrid stresses. Consider a fluctuation defined by $\bm{x}' = \bm{x} - \overline{\bm{x}}$. Taking the derivative and substituting
\begin{align}
    \bm{x}'_t = \bm{x}_t - \overline{\bm{x}}_t &= \bm{f}(\bm{x}) - \bm{f}(\overline{\bm{x}}) - \bm{\tau}(\overline{\bm{x}},\bm{x}) \\
    &= \bm{f}(\overline{\bm{x}} + \bm{x}') - \overline{\bm{f}(\overline{\bm{x}}+ \bm{x}')}
\end{align}
Expand $\bm{f}(\bm{x})$ for $\bm{x} \approx \overline{\bm{x}}$ using a Taylor series
\begin{align}
    \bm{f}(\overline{\bm{x}} + \bm{x}') &= \bm{f}(\overline{\bm{x}}) + \bm{J}(\overline{\bm{x}}) \cdot \bm{x}' + \frac{1}{2} \cdot \bm{x}'^T \cdot \bm{H}(\overline{\bm{x}}) \cdot \bm{x}'
\end{align}
where $\bm{J}$ and $\bm{H}$ are the Jacobian and Hessian matrices. Applying the filtering operator yields
\begin{align}
    \overline{\bm{f}(\overline{\bm{x}} + \bm{x}')} &= \overline{\bm{f}(\overline{\bm{x}}) + \bm{J}(\overline{\bm{x}}) \cdot \bm{x}' + \frac{1}{2} \cdot \bm{x}'^T \cdot \bm{H}(\overline{\bm{x}}) \cdot \bm{x}'} \\
    \overline{\bm{f}(\bm{x})} &= \bm{f}(\overline{\bm{x}}) + \overline{\frac{1}{2} \cdot \bm{x}'^T \cdot \bm{H}(\overline{\bm{x}}) \cdot \bm{x}'}
\end{align}
So, for the subgrid stresses, we have
\begin{align}
    \bm{\tau}(\overline{\bm{x}}, \bm{x}) =  \overline{\frac{1}{2} \cdot (\bm{x} - \overline{\bm{x}})^T \cdot \bm{H}(\overline{\bm{x}}) \cdot (\bm{x} - \overline{\bm{x}})}
\end{align}
Assuming linear dynamics for the fluctuations $\bm{x}'$ and small filter width $\Delta$ gives
\begin{align}
    &\bm{x}'_t = \bm{J}(\overline{\bm{x}}) \cdot \bm{x}' \\
    &\bm{x}' = e^{\bm{J}\Delta} \bm{x}'_0 \approx (\bm{I} + \bm{J}\Delta)\bm{x}'_0
\end{align}
Substitute and apply the filter operator
\begin{align}
    \frac{1}{2} \cdot \bm{x}'^T \cdot \bm{H}(\overline{\bm{x}}) \cdot \bm{x}' &= \frac{1}{2} [(\bm{I} + \bm{J}\Delta)\bm{x}'_0]^T \cdot \bm{H} \cdot (\bm{I} + \bm{J}\Delta)\bm{x}'_0 \\
    \frac{1}{2} \cdot \bm{x}'^T \cdot \bm{H}(\overline{\bm{x}}) \cdot \bm{x}' &= \frac{1}{2} \left( \bm{x}'^T_0 \bm{H} \bm{x}'_0 + \bm{x}'^T_0 (\bm{J}^T \bm{H} + \bm{H}\bm{J}) \bm{x}'_0 \lambda + \bm{x}'^T_0 \bm{J}^T \bm{H} \bm{J} \bm{x}'_0 \lambda^2 \right) \\
    \overline{\frac{1}{2} \cdot \bm{x}'^T \cdot \bm{H}(\overline{\bm{x}}) \cdot \bm{x}'} &= \frac{1}{2\Delta} \int_{-\Delta/2}^{\Delta/2} \left(\bm{x}'^T_0 \bm{H} \bm{x}'_0 + \bm{x}'^T_0 \left(\bm{J}^T \bm{H} + \bm{H} \bm{J}\right)\bm{x}'_0 \tau + \bm{x}'^T_0  \bm{J}^T \bm{H} \bm{J} \bm{x}'_0 \tau^2 \right) d\tau \\ 
    &= \frac{1}{2} \bm{x}'^T_0 \left[ \bm{H} + \frac{\Delta^2}{12} \bm{J}^T \bm{H} \bm{J} \right] \bm{x'}_0
\end{align}

\newpage
\printbibliography

@article{volkmann2024scalable,
  title={A scalable generative model for dynamical system reconstruction from neuroimaging data},
  author={Volkmann, E. and Br{\"a}ndle, A. and Durstewitz, D. and Koppe, G.},
  journal={Advances in Neural Information Processing Systems},
  volume={37},
  pages={80328--80362},
  year={2024}
}

@misc{flowsanddiffusions2025,
    author       = {Holderrieth, P. and Erives, E.},
    title        = {Introduction to Flow Matching and Diffusion Models},
    year         = {2025},
    url          = {https://diffusion.csail.mit.edu/}
  }

@misc{ramdas2015wassersteinsampletestingrelated,
      title={On Wasserstein Two Sample Testing and Related Families of Nonparametric Tests}, 
      author={Ramdas, A. and Garcia, N. and Cuturi, M.},
      year={2015},
      eprint={1509.02237},
      archivePrefix={arXiv},
      primaryClass={math.ST},
      url={https://arxiv.org/abs/1509.02237}, 
}

@article{kolouri2017optimalmass,
  author={Kolouri, S. and Park, S. R. and Thorpe, M. and Slepcev, D. and Rohde, G. K.},
  journal={IEEE Signal Processing Magazine}, 
  title={Optimal Mass Transport: Signal processing and machine-learning applications}, 
  year={2017},
  volume={34},
  number={4},
  pages={43-59},
  doi={10.1109/MSP.2017.2695801}
}

@misc{ho2022classifierfreediffusionguidance,
      title={Classifier-Free Diffusion Guidance}, 
      author={Ho, J. and Salimans, T.},
      year={2022},
      eprint={2207.12598},
      archivePrefix={arXiv},
      primaryClass={cs.LG},
      url={https://arxiv.org/abs/2207.12598}, 
}

@article{jacobsen2023cocogen,
  title={Cocogen: Physically-consistent and conditioned score-based generative models for forward and inverse problems},
  author={Jacobsen, C. and Zhuang, Y. and Duraisamy, K.},
  journal={arXiv preprint arXiv:2312.10527},
  year={2023}
}

@article{song2020score,
  title={Score-based generative modeling through stochastic differential equations},
  author={Song, Y. and Sohl-Dickstein, J. and Kingma, D. P. and Kumar, A. and Ermon, S. and Poole, B.},
  journal={arXiv preprint arXiv:2011.13456},
  year={2020}
}

@article{ho2020denoising,
  title={Denoising diffusion probabilistic models},
  author={Ho, J. and Jain, A. and Abbeel, P.},
  journal={Advances in neural information processing systems},
  volume={33},
  pages={6840--6851},
  year={2020}
}

@article{finn2024generative,
author = {Finn, T. B. and Durand, C. and Farchi, A. and Bocquet, M. and Rampal, P. and Carrassi, A.},
title = {Generative Diffusion for Regional Surrogate Models From Sea-Ice Simulations},
journal = {Journal of Advances in Modeling Earth Systems},
volume = {16},
number = {10},
doi = {https://doi.org/10.1029/2024MS004395},
year = {2024}
}

@article{adrian_stochastic_1990,
  title = {Stochastic {{Estimation}} of {{Sub-Grid Scale Motions}}},
  author = {Adrian, R. J.},
  year = {1990},
  month = may,
  journal = {Applied Mechanics Reviews},
  volume = {43},
  number = {5S},
  pages = {S214-218},
  issn = {0003-6900, 2379-0407},
  doi = {10.1115/1.3120809},
  langid = {english}
}

@article{arnold_stochastic_2013,
  title = {Stochastic Parametrizations and Model Uncertainty in the {{Lorenz}} '96 System},
  author = {Arnold, H. M. and Moroz, I. M. and Palmer, T. N.},
  year = {2013},
  month = may,
  journal = {Philosophical Transactions of the Royal Society A: Mathematical, Physical and Engineering Sciences},
  volume = {371},
  number = {1991},
  pages = {20110479},
  publisher = {Royal Society},
  doi = {10.1098/rsta.2011.0479}
}

@article{bayram_numerical_2018,
  title = {Numerical Methods for Simulation of Stochastic Differential Equations},
  author = {Bayram, M. and Partal, T. and Orucova Buyukoz, G.},
  year = {2018},
  month = jan,
  journal = {Advances in Difference Equations},
  volume = {2018},
  number = {1},
  pages = {17},
  issn = {1687-1847},
  doi = {10.1186/s13662-018-1466-5}
}

@article{bose_wallmodeled_2018,
  title = {Wall-{{Modeled Large-Eddy Simulation}} for {{Complex Turbulent Flows}}},
  author = {Bose, S. T. and Park, G. I.},
  year = {2018},
  month = jan,
  journal = {Annual Review of Fluid Mechanics},
  volume = {50},
  number = {Volume 50, 2018},
  pages = {535--561},
  publisher = {Annual Reviews},
  issn = {0066-4189, 1545-4479},
  doi = {10.1146/annurev-fluid-122316-045241},
  langid = {english}
}

@article{bodart_wallmodeled_2011,
  title = {Wall-Modeled Large Eddy Simulation in Complex Geometries with Application to High-Lift Devices},
  author = {Bodart, J. and Larsson, J.},
  year = {2011},
  langid = {english}
}

@article{chorin_optimal_2002,
  title = {Optimal Prediction with Memory},
  author = {Chorin, A. J. and Hald, O. H. and Kupferman, R.},
  year = {2002},
  month = jun,
  journal = {Physica D: Nonlinear Phenomena},
  volume = {166},
  number = {3-4},
  pages = {239--257},
  issn = {01672789},
  doi = {10.1016/S0167-2789(02)00446-3},
  copyright = {https://www.elsevier.com/tdm/userlicense/1.0/},
  langid = {english}
}

@article{christensen_simulating_2015,
  title = {Simulating Weather Regimes: Impact of Stochastic and Perturbed Parameter Schemes in a Simple Atmospheric Model},
  shorttitle = {Simulating Weather Regimes},
  author = {Christensen, H. M. and Moroz, I. M. and Palmer, T. N.},
  year = {2015},
  month = apr,
  journal = {Climate Dynamics},
  volume = {44},
  number = {7},
  pages = {2195--2214},
  issn = {1432-0894},
  doi = {10.1007/s00382-014-2239-9},
  langid = {english}
}

@article{cruzeiro_stochastic_2020,
  title = {Stochastic {{Approaches}} to {{Deterministic Fluid Dynamics}}: {{A Selective Review}}},
  shorttitle = {Stochastic {{Approaches}} to {{Deterministic Fluid Dynamics}}},
  author = {Cruzeiro, A. B.},
  year = {2020},
  month = mar,
  journal = {Water},
  volume = {12},
  number = {3},
  pages = {864},
  issn = {2073-4441},
  doi = {10.3390/w12030864},
  copyright = {https://creativecommons.org/licenses/by/4.0/},
  langid = {english}
}

@article{dietrich_learning_2023,
  title = {Learning Effective Stochastic Differential Equations from Microscopic Simulations: {{Linking}} Stochastic Numerics to Deep Learning},
  shorttitle = {Learning Effective Stochastic Differential Equations from Microscopic Simulations},
  author = {Dietrich, F. and Makeev, A. and Kevrekidis, G. and Evangelou, N. and Bertalan, T. and Reich, S. and Kevrekidis, I. G.},
  year = {2023},
  month = feb,
  journal = {Chaos: An Interdisciplinary Journal of Nonlinear Science},
  volume = {33},
  number = {2},
  pages = {023121},
  issn = {1054-1500},
  doi = {10.1063/5.0113632}
}

@article{garnier_large_2002,
  title = {Large {{Eddy Simulation}} of {{Shock}}/{{Boundary-Layer Interaction}}},
  author = {Garnier, E. and Sagaut, P. and Deville, M.},
  year = {2002},
  month = oct,
  journal = {AIAA Journal},
  volume = {40},
  number = {10},
  pages = {1935--1944},
  publisher = {{American Institute of Aeronautics and Astronautics}},
  issn = {0001-1452},
  doi = {10.2514/2.1552},
}

@book{garnier_large_2009,
  title = {Large {{Eddy Simulation}} for {{Compressible Flows}}},
  author = {Garnier, E. and Adams, N. and Sagaut, P.},
  editor = {Chattot, J.-J. and Colella, P. and Eist, W. and Glowinski, R. and Hussaini, Y. and Joly, P. and Keller, H. B. and Marsden, J. E. and Meiron, D. I. and Pironneau, O. and Quarteroni, A. and Rappaz, J. and Rosner, R. and Sagaut, P. and Seinfeld, J. H. and Szepessy, A. and Wheeler, M. F.},
  year = {2009},
  series = {Scientific {{Computation}}},
  publisher = {Springer Netherlands},
  address = {Dordrecht},
  doi = {10.1007/978-90-481-2819-8},
  langid = {english}
}

@article{germano_turbulence_1992,
  title = {Turbulence: The Filtering Approach},
  shorttitle = {Turbulence},
  author = {Germano, M.},
  year = {1992},
  month = may,
  journal = {Journal of Fluid Mechanics},
  volume = {238},
  pages = {325--336},
  issn = {0022-1120, 1469-7645},
  doi = {10.1017/S0022112092001733},
  langid = {english}
}

@article{givon_extracting_2004,
  title = {Extracting Macroscopic Dynamics: Model Problems and Algorithms},
  shorttitle = {Extracting Macroscopic Dynamics},
  author = {Givon, D. and Kupferman, R. and Stuart, A.},
  year = {2004},
  month = nov,
  journal = {Nonlinearity},
  volume = {17},
  number = {6},
  pages = {R55-R127},
  issn = {0951-7715, 1361-6544},
  doi = {10.1088/0951-7715/17/6/R01}
}

@article{gottwald_datadriven_2016,
  title = {A Data-Driven Method for the Stochastic Parametrisation of Subgrid-Scale Tropical Convective Area Fraction},
  author = {Gottwald, G. A. and Peters, K. and Davies, L.},
  year = {2016},
  journal = {Quarterly Journal of the Royal Meteorological Society},
  volume = {142},
  number = {694},
  pages = {349--359},
  issn = {1477-870X},
  doi = {10.1002/qj.2655},
  langid = {english}
}

@incollection{gottwald_stochastic_2017,
  title = {Stochastic {{Climate Theory}}},
  booktitle = {Nonlinear and {{Stochastic Climate Dynamics}}},
  author = {Gottwald, G. A. and Crommelin, D. T. and Franzke, C. L. E.},
  editor = {Franzke, Christian L. E. and O'Kane, Terence J.},
  year = {2017},
  pages = {209--240},
  publisher = {Cambridge University Press},
  address = {Cambridge},
  doi = {10.1017/9781316339251.009}
}

@article{grudzien_numerical_2020,
  title = {On the Numerical Integration of the {{Lorenz-96}} Model, with Scalar Additive Noise, for Benchmark Twin Experiments},
  author = {Grudzien, C. and Bocquet, M. and Carrassi, A.},
  year = {2020},
  month = apr,
  journal = {Geoscientific Model Development},
  volume = {13},
  number = {4},
  pages = {1903--1924},
  publisher = {Copernicus GmbH},
  issn = {1991-959X},
  doi = {10.5194/gmd-13-1903-2020},
  langid = {english}
}

@article{hasselmann_stochastic_1976,
  title = {Stochastic Climate Models {{Part I}}. {{Theory}}},
  author = {Hasselmann, K.},
  year = {1976},
  journal = {Tellus},
  volume = {28},
  number = {6},
  pages = {473--485},
  issn = {2153-3490},
  doi = {10.1111/j.2153-3490.1976.tb00696.x},
  copyright = {1976 Blackwell Munksgaard},
  langid = {english}
}

@article{higham_algorithmic_2001,
  title = {An {{Algorithmic Introduction}} to {{Numerical Simulation}} of {{Stochastic Differential Equations}}},
  author = {Higham, D. J.},
  year = {2001},
  month = jan,
  journal = {SIAM Review},
  volume = {43},
  number = {3},
  pages = {525--546},
  issn = {0036-1445, 1095-7200},
  doi = {10.1137/S0036144500378302},
  langid = {english}
}

@book{khasminskii_stochastic_2012,
  title = {Stochastic {{Stability}} of {{Differential Equations}}},
  author = {Khasminskii, R.},
  year = {2012},
  series = {Stochastic {{Modelling}} and {{Applied Probability}}},
  volume = {66},
  publisher = {Springer Berlin Heidelberg},
  address = {Berlin, Heidelberg},
  doi = {10.1007/978-3-642-23280-0},
  langid = {english}
}

@book{kloeden_numerical_1992,
  title = {Numerical {{Solution}} of {{Stochastic Differential Equations}}},
  author = {Kloeden, P. E. and Platen, E.},
  year = {1992},
  publisher = {Springer Berlin Heidelberg},
  address = {Berlin, Heidelberg},
  doi = {10.1007/978-3-662-12616-5},
  langid = {english}
}

@article{leith_stochastic_1990,
  title = {Stochastic Backscatter in a Subgrid-Scale Model: {{Plane}} Shear Mixing Layer},
  shorttitle = {Stochastic Backscatter in a Subgrid-Scale Model},
  author = {Leith, C. E.},
  year = {1990},
  month = mar,
  journal = {Physics of Fluids A},
  volume = {2},
  pages = {297--299},
  issn = {0899-82131070-6631},
  doi = {10.1063/1.857779}
}

@article{lilly_numerical_1962,
  title = {On the Numerical Simulation of Buoyant Convection},
  author = {Lilly, D. K.},
  year = {1962},
  month = may,
  journal = {Tellus},
  volume = {14},
  number = {2},
  pages = {148--172},
  issn = {00402826, 21533490},
  doi = {10.1111/j.2153-3490.1962.tb00128.x},
  langid = {english}
}

@book{liu_stochastic_2019,
  title = {Stochastic {{Stability}} of {{Differential Equations}} in {{Abstract Spaces}}},
  author = {Liu, K.},
  year = {2019},
  series = {London {{Mathematical Society Lecture Note Series}}},
  publisher = {Cambridge University Press},
  address = {Cambridge},
  doi = {10.1017/9781108653039}
}

@article{majda_mathematical_2001,
  title = {A Mathematical Framework for Stochastic Climate Models},
  author = {Majda, A. J. and Timofeyev, I. and Vanden Eijnden, E.},
  year = {2001},
  month = aug,
  journal = {Communications on Pure and Applied Mathematics},
  volume = {54},
  number = {8},
  pages = {891--974},
  issn = {0010-3640, 1097-0312},
  doi = {10.1002/cpa.1014},
  langid = {english}
}

@article{marstorp_stochastic_2007,
  title = {A Stochastic Subgrid Model with Application to Turbulent Flow and Scalar Mixing},
  author = {Marstorp, L. and Brethouwer, G. and Johansson, A. V.},
  year = {2007},
  month = mar,
  journal = {Physics of Fluids},
  volume = {19},
  number = {3},
  pages = {035107},
  issn = {1070-6631, 1089-7666},
  doi = {10.1063/1.2711477},
  langid = {english}
}

@article{mori_transport_1965,
  title = {Transport, {{Collective Motion}}, and {{Brownian Motion}}},
  author = {Mori, H.},
  year = {1965},
  month = mar,
  journal = {Progress of Theoretical Physics},
  volume = {33},
  number = {3},
  pages = {423--455},
  issn = {0033-068X},
  doi = {10.1143/PTP.33.423},
  langid = {english}
}

@article{moser_statistical_2021,
  title = {Statistical {{Properties}} of {{Subgrid-Scale Turbulence Models}}},
  author = {Moser, R. D. and Haering, S. W. and Yalla, G. R.},
  year = {2021},
  month = jan,
  journal = {Annual Review of Fluid Mechanics},
  volume = {53},
  number = {Volume 53, 2021},
  pages = {255--286},
  publisher = {Annual Reviews},
  issn = {0066-4189, 1545-4479},
  doi = {10.1146/annurev-fluid-060420-023735},
  langid = {english}
}

@book{ottinger_stochastic_1996,
  title = {Stochastic {{Processes}} in {{Polymeric Fluids}}: {{Tools}} and {{Examples}} for {{Developing Simulation Algorithms}}},
  shorttitle = {Stochastic {{Processes}} in {{Polymeric Fluids}}},
  author = {{\"O}ttinger, H. C.},
  year = {1996},
  publisher = {Springer Berlin Heidelberg},
  address = {Berlin, Heidelberg},
  doi = {10.1007/978-3-642-58290-5},
  copyright = {https://www.springernature.com/gp/researchers/text-and-data-mining},
  langid = {english}
}

@article{palmer_stochastic_2019,
  title = {Stochastic Weather and Climate Models},
  author = {Palmer, T. N.},
  year = {2019},
  month = may,
  journal = {Nature Reviews Physics},
  volume = {1},
  number = {7},
  pages = {463--471},
  issn = {2522-5820},
  doi = {10.1038/s42254-019-0062-2},
  langid = {english}
}

@book{paul_stochastic_2013,
  title = {Stochastic {{Processes}}: {{From Physics}} to {{Finance}}},
  shorttitle = {Stochastic {{Processes}}},
  author = {Paul, W. and Baschnagel, J.},
  year = {2013},
  publisher = {Springer International Publishing},
  address = {Heidelberg},
  doi = {10.1007/978-3-319-00327-6},
  copyright = {https://www.springernature.com/gp/researchers/text-and-data-mining},
  langid = {english}
}

@article{perezhogin_generative_2023,
  title = {Generative {{Data-Driven Approaches}} for {{Stochastic Subgrid Parameterizations}} in an {{Idealized Ocean Model}}},
  author = {Perezhogin, P. and Zanna, L. and {Fernandez-Granda}, C.},
  year = {2023},
  journal = {Journal of Advances in Modeling Earth Systems},
  volume = {15},
  number = {10},
  pages = {e2023MS003681},
  issn = {1942-2466},
  doi = {10.1029/2023MS003681},
  copyright = {{\copyright} 2023 The Authors. Journal of Advances in Modeling Earth Systems published by Wiley Periodicals LLC on behalf of American Geophysical Union.},
  langid = {english}
}

@article{rasam_stochastic_2014,
  title = {A Stochastic Extension of the Explicit Algebraic Subgrid-Scale Models},
  author = {Rasam, A. and Brethouwer, G. and Johansson, A. V.},
  year = {2014},
  month = may,
  journal = {Physics of Fluids},
  volume = {26},
  number = {5},
  pages = {055113},
  issn = {1070-6631, 1089-7666},
  doi = {10.1063/1.4879436},
  langid = {english}
}

@book{sagaut_large_2006,
  title = {Large Eddy Simulation for Incompressible Flows: An Introduction},
  shorttitle = {Large Eddy Simulation for Incompressible Flows},
  author = {Sagaut, P.},
  year = {2006},
  series = {Scientific Computation},
  edition = {3rd ed},
  publisher = {Springer},
  address = {Berlin ; New York},
  langid = {english},
  lccn = {MLCM 2006/40205 (T)}
}

@article{sagaut_large_2009,
  title = {Large Eddy Simulation for Aerodynamics: Status and Perspectives},
  shorttitle = {Large Eddy Simulation for Aerodynamics},
  author = {Sagaut, P. and Deck, S.},
  year = {2009},
  month = jul,
  journal = {Philosophical Transactions of the Royal Society A: Mathematical, Physical and Engineering Sciences},
  volume = {367},
  number = {1899},
  pages = {2849--2860},
  publisher = {Royal Society},
  doi = {10.1098/rsta.2008.0269}
}

@article{schumann_stochastic_1997,
  title = {Stochastic Backscatter of Turbulence Energy and Scalar Variance by Random Subgrid-Scale Fluxes},
  author = {Schumann, U. and Launder, B. E.},
  year = {1997},
  month = jan,
  journal = {Proceedings of the Royal Society of London. Series A: Mathematical and Physical Sciences},
  volume = {451},
  number = {1941},
  pages = {293--318},
  publisher = {Royal Society},
  doi = {10.1098/rspa.1995.0126}
}

@techreport{slotnick_cfd_2014,
  title = {{{CFD Vision}} 2030 {{Study}}: {{A Path}} to {{Revolutionary Computational Aerosciences}}},
  shorttitle = {{{CFD Vision}} 2030 {{Study}}},
  author = {Slotnick, J. P. and Khodadoust, A. and Alonso, J. and Darmofal, D. and Gropp, W. and Lurie, E. and Mavriplis, D. J.},
  year = {2014},
  month = mar,
  number = {NF1676L-18332}
}

@article{smagorinsky_general_1963,
  title = {{{GENERAL CIRCULATION EXPERIMENTS WITH THE PRIMITIVE EQUATIONS}}: {{I}}. {{THE BASIC EXPERIMENT}}*},
  shorttitle = {{{GENERAL CIRCULATION EXPERIMENTS WITH THE PRIMITIVE EQUATIONS}}},
  author = {Smagorinsky, J.},
  year = {1963},
  month = mar,
  journal = {Monthly Weather Review},
  volume = {91},
  number = {3},
  pages = {99--164},
  issn = {0027-0644, 1520-0493},
  doi = {10.1175/1520-0493(1963)091<0099:GCEWTP>2.3.CO;2},
  langid = {english}
}

@article{vreman_eddyviscosity_2004,
  title = {An Eddy-Viscosity Subgrid-Scale Model for Turbulent Shear Flow: {{Algebraic}} Theory and Applications},
  shorttitle = {An Eddy-Viscosity Subgrid-Scale Model for Turbulent Shear Flow},
  author = {Vreman, A. W.},
  year = {2004},
  month = oct,
  journal = {Physics of Fluids},
  volume = {16},
  number = {10},
  pages = {3670--3681},
  issn = {1070-6631, 1089-7666},
  doi = {10.1063/1.1785131},
  langid = {english}
}

@article{deardorff_numerical_1970,
  title = {A Numerical Study of Three-Dimensional Turbulent Channel Flow at Large {{Reynolds}} Numbers},
  author = {Deardorff, J. W.},
  year = {1970},
  month = apr,
  journal = {Journal of Fluid Mechanics},
  volume = {41},
  number = {2},
  pages = {453--480},
  issn = {0022-1120, 1469-7645},
  doi = {10.1017/S0022112070000691},
  copyright = {https://www.cambridge.org/core/terms},
  langid = {english}
}

@article{zwanzig_memory_1961,
  title = {Memory {{Effects}} in {{Irreversible Thermodynamics}}},
  author = {Zwanzig, R.},
  year = {1961},
  month = nov,
  journal = {Physical Review},
  volume = {124},
  number = {4},
  pages = {983--992},
  publisher = {American Physical Society},
  doi = {10.1103/PhysRev.124.983}
}

\end{document}